\documentclass[10pt]{article}
\usepackage{bbm}
\usepackage{amsmath}
\usepackage{amsthm}
\usepackage{amssymb} 
\usepackage{mathrsfs}
\usepackage{graphicx}
\usepackage{authblk}
\usepackage[backref=none,hypertexnames=false,colorlinks=true,urlcolor=blue,linkcolor=blue,citecolor=blue]{hyperref}
\usepackage[numbers,comma,square,sort&compress]{natbib}
\usepackage[letterpaper,text={6.5in,9in},centering]{geometry}

\graphicspath{{Figs}}

\makeatletter\@addtoreset{equation}{section}\makeatother

\def\Ab{\mathbf A}
\def\Ib{\mathbf I}
\def\Qb{\mathbf Q}

\def\xb{\mathbf x}

\def\wb{\mathbf w}
\def\pb{\mathbf p}
\def\qb{\mathbf q}
\def\yb{\mathbf y}
\def\zb{\mathbf z}

\newcommand{\Pcov}{\mathbf{P}}

 {\begin{trivlist}\item[]\textbf{Acknowledgments.}}{\end{trivlist}}

\begin{document}
\title{Influence network reconstruction from discrete time-series of count data modelled by multidimensional Hawkes processes}
\author[1]{Naratip Santitissadeekorn}
\author[2]{Martin Short}
\author[1,3]{David J. B. Lloyd}

\affil[1]{\small Department of Mathematics, University of Surrey, Guildford, UK}
\affil[2]{\small Department of Mathematics, Georgia Institute of Technology, USA}
\affil[3]{\small Centre for Criminology, University of Surrey, Guildford, UK}
\date{\today}
\maketitle

\begin{abstract}
Identifying key influencers from time series data without a known prior network structure is a challenging problem in various applications, from crime analysis to social media. While much work has focused on event-based time series (timestamp) data, fewer methods address count data, where event counts are recorded in fixed intervals. We develop network inference methods for both batched and sequential count data. Here the strong network connection represents the key influences among the nodes. We introduce an ensemble-based algorithm, rooted in the expectation-maximization (EM) framework, and demonstrate its utility to identify node dynamics and connections through a discrete-time Cox or Hawkes process. For the linear multidimensional Hawkes model, we employ a minimization-majorization (MM) approach, allowing for parallelized inference of networks.
For sequential inference, we use a second-order approximation of the Bayesian inference problem. Under certain assumptions, a rank-1 update for the covariance matrix reduces computational costs. We validate our methods on synthetic data and real-world datasets, including email communications within European academic communities. Our approach effectively reconstructs underlying networks, accounting for both excitation and diffusion influences. This work advances network reconstruction from count data in real-world scenarios.

\end{abstract}

\noindent
{\bf Keywords:}
  Influence network, causal network, nonlinear filtering, Hawkes process, count data, expectation maximization

\noindent
{\bf AMS Classification:}
  62F15, 62F30, 62M20

\section{Introduction}

This work is motivated by conventional applications of continuous-time Hawkes processes utilized to model the temporal clustering as well as mutual excitation network driven by the timestamp data, i.e. the times of events. The continuous-time point-process Hawkes model was first introduced by Hawkes \cite{Hawkes} to capture a self-excitation process, used particularly for seismic events~\cite{Ogata1999SeismicityAT}. Since then, it has been extended to multivariate Hawkes process to model the mutual excitation structure or influence network structure. This development has led to emerging applications of point-process Hawkes model to seismic analysis~\cite{Veen08,Ogata1999SeismicityAT}, urban crime analysis~\cite{mohler2012geographic, Mohler13, short2008, YuanSIAMdatasci} social network analysis~\cite{Lima23, Zhou13a,ZhaoEHRL15,HsinYU17,Kobayashi16,Chen18,Eric16,zipkin_schoenberg_coronges_bertozzi_2016}, financial time-series analysis~\cite{Bacry12,Bacry15,Heidar23,Hawkes18finance}, contagious disease network~\cite{Unwin21,Choi2015ConstructingDN,Diouane} and deep learning network~\cite{Deeplearn_Shang_Sun_2019,eqmarginalLH,Deeplearn_wangg16,Deeplearn_Xiao,Deeplearn_zuo20a}.
\par
To reconstruct the influence network from a time series of timestamp data, most of the above work typically used the Expectation-Maximization (EM) or Minimization-Majorization (MM) framework to construct a surrogate function (i.e., tight upper bound function) for the negative log-likelihood function. 
The main advantage of this approach is that it may help decouple the parameter space when optimizing the surrogate function, accelerating the computational task. For a simple excitation kernel, such as the exponential decay kernel, a closed-form method for the parameter update can be derived. A nonparametric excitation kernel can also be used within the EM and MM approach, where the Euler-Lagrange equation can be derived for the optimization of the surrogate function ~\cite{LewisMohler} and the regularization to promote sparsity~\cite{Zhao15}. Other techniques have also been developed to estimate non-parametric kernels; see for instance~\cite{Bacry12, Lemonnier14, Eichler16}. 
A fully-Bayesian, parallel inference algorithm was also developed in \cite{Linderman14} to model the excitation structure by random graph models which allows conjugate prior for efficient inference via Markov chain Monte Carlo.
\par
The influence network within multi-dimensional Hawkes models can also be linked to Granger-causality in temporal point processes \cite{Granger96}. In the context of this framework, an event generated by $x_j$ is considered to ``Granger-cause" the event associated with $x_i$ if the likelihood function of events in $x_i(t)$, given the history of all events up to time $t$, decreases when the history of events generated by $x_j$ is omitted. The application of the multivariate Hawkes model in the context of Granger causality provides interpretability of the results. It was demonstrated in~\cite{Eichler16} that $x_j$ does not Granger-cause $x_i$ when the (pairwise) excitation kernel used by $x_j$ to 'excite' $x_i$ is zero. Applications of Hawkes model to discover Granger-causality were investigated with real-world data in~\cite{Xu16, Bacry17,Kimetal11}. However, this work is limited to certain conditional independence of the excitation process, while some data may exhibit inhibition or interaction. Additionally, to the best of our knowledge, the connection between the excitation (or influence) network and Granger causality has not been extended to the case of count data modeled by the discrete-time Hawkes process. Therefore, the influence network derived from count data may or may not correspond to the causal network.
\par
Similar to the timestamp data, a time series of count data may exhibit self-excitation, wherein a high count is often followed by several higher counts, e.g., in a time series of epileptic seizure counts~\cite{Albert91}. This data is often easier to collect and more common in applications, for instance in epidemiology, where one can only sensibly collate count data. Moreover, multiple time series may demonstrate an ``influencing" characteristic, where a high count from one series is followed by high counts in others. For instance, a cluster of earthquakes in a particular region could trigger seismic activities in adjacent regions, while incidents occurring in one area of a city could lead to similar occurrences in other areas, e.g., urban crimes~\cite{short2008,short2010}. We can conceptualize the sources of these multiple time series as nodes in a network, and by uncovering the influence structure of such a network, we can gain insights into the evolving dynamics of the network over time, such as the emergence of synchronisation of node dynamics. However, there are no methods (to the best of the author's knowledge) for carrying out the network reconstruction problem from time series count data. 
\par
The primary aim of this work is to identify an influence network from a time series consisting of count data that opens up more real-world applications of network reconstruction via Hawkes-type processes. The count data inference problem is significantly more challenging since data is aggregated and therefore there is a loss of information relative to the time-stamp data which needs to be accounted for in the inference problem. Recent research \cite{Rizoiu22} has developed a novel Hawkes-type non-homogeneous Poisson process for count data, with a deterministic intensity function equal to the expectation of a Hawkes process intensity and parameters are fitted based on the interval-censored likelihood function derived therein. Here we adopt different models that utilize a discrete-time, multidimensional Cox or Hawkes process to model the excitation effects among nodes driven by count data. Within this framework, the magnitude of the influence can be quantified by the excitation rate parameters incorporated into the model. Therefore, the task of identifying the influence network reduces to estimating the parameters of the Cox or Hawkes model. This work presents three distinct methodologies for parameter estimation for count data problems: (1) Ensemble-based EM, (2) Minimization-Majorization (MM) technique and (3) a sequential algorithm based on approximate second-order filtering.

For the ensemble-based EM, the Hawkes process can take a general ``state-space" form (e.g., doubly stochastic Poisson point process) that consists of state dynamical system and observation equation. Within the context of EM, the ``missing" data is the unobserved sample path of the state. Therefore, the E-step requires a Monte Carlo sampling of the sample paths. This approach was previously used in the context of the model identification for the Kalman filter \cite{shumwaystoffer82,Godsill2004MonteCS}. We demonstrate how this idea can be applied to network reconstruction of a small network. 

The MM algorithm is developed for batch data inference. For this algorithm, we limit the node dynamic model to a discrete-time dynamical system analogous to the exponential-decay kernel of the multivariate Hawkes process. We show how one can derive an iterative method that minimises a surrogate function such that the surrogate function is a ``tight" upper bound of the negative log-likelihood function for count data. We show how the iterative method can be parallelised. 

The sequential algorithm based on approximate second-order filtering, called the extended Poisson-Kalman filter (ExPKF), is derived by approximating the mean and covariance matrix of the posterior density for the same dynamical systems Hawkes model used for the MM algorithm. We show how this leads to an efficient method under an assumption, where one can use a rank-1 update for the update of the covariance matrix with parallelisation.

Our main contributions are the development of the foundations for a systematic approach to dealing with network influence reconstruction from count data. We present methods that deal with either complex state-space models for small networks or linear Hawkes processes for large networks. The ensemble-based EM method also captures uncertainty quantification of the intensity estimate, while ExPKF provides a second-order moment for the network estimate. Uncertainty quantification is very important in network reconstruction applications so that one can gain an understanding of the uncertainty of a link between two nodes occurring. This work opens up new avenues of research involving count data collected on networks, and the development of new methods for more general or other types of stochastic processes. 

The paper is outlined as follows. In Section \ref{sec:Ensemble-based EM}, we develop the ensemble-based EM algorithm for small networks. We then focus on a moderate-size network for a count-data model, inspired by the discretization of the exponential decay kernel of the continuous-time Hawkes process. The MM algorithm is developed for batch data inference in Section~\ref{sec:MM} and the extended Poisson-Kalman filter (ExPKF) is derived in Section~\ref{s:seq_DA}. In section~\ref{s:results}, we demonstrate the validity of the proposed methods to reconstruct the influence network with various numerical experiments with known ground truths. We also demonstrate the utility of the method on large real-world email network data in section~\ref{s:email} and conclude in section~\ref{s:con}.

\section{Ensemble-based EM}\label{sec:Ensemble-based EM}
We are interested in an inhomogeneous Poisson point process on a network with $m$ nodes, where the conditional intensity $\lambda_k^i$ at the $i$-th node is assumed to be a constant in the $k-$th time interval $(t_{k},t_{k+1})$. In other words, if $\Delta N_k^i$ is the number of events observed for the $i-$th node at the $k-$th time interval, we assume that $Pr(\Delta N_k^i\mid\lambda_k^i)$ is a Poisson probability with mean $\lambda_k^i\delta t_k$ where $\delta t_k=t_{k+1}-t_{k}$.  The intensity $\lambda_k^i$ depends on a $n$-dimensional parameter vector, $\theta^i:=[\theta^{i,1},\ldots,\theta^{i,n}]^\top$.
We concatenate all vectors $\theta^i$ to form a parameter vector $\theta$, i.e., $\theta:=[\theta^{1};\cdots;\theta^{m}]$.
\par
For simplicity, all the time intervals are assumed to have the same length $\delta t$. At any given time step $k$, we assume conditional independence so that the conditional joint density is given by 
\begin{equation}\label{eq:PoiLikelihood}
p(\underbrace{\Delta N_{k}^1,\ldots,\Delta N_{k}^m}_{\equiv\Delta N_k}|\lambda^1_k,\cdots,\lambda^m_k)\propto\prod_{i=1}^m(\lambda^i_k)^{\Delta N^i_k}\exp(-\lambda^i_k\delta t).
\end{equation}
Let $\Delta N_{1:K}:=\left[\Delta N_1,\ldots,\Delta N_K\right]$ denote time-series of count data up to the time step $K$ for all nodes. The log-likelihood function is then given by
\begin{equation}\label{eq:LL}
    \mathbf{L}(\theta):=\log p(\Delta N_{1:K}\mid\theta)=\sum_{i=1}^m\sum_{k=1}^K\log(\lambda^i_k(\theta^i))\Delta N_k^i-\delta t\sum_{i=1}^m\sum_{k=1}^K\lambda^i_k(\theta^i)+\mathcal{C},
\end{equation}
where $\mathcal{C}$ is independent of $\theta$.
The maximum likelihood method estimates the model parameter vector $\theta$ (in a parameter space $\Theta$) by maximizing the log-likelihood function 
\begin{equation}
\widehat{\theta}:=\arg\max_{\theta\in\Theta}~\mathbf{L}(\theta).
\end{equation}
\par
We assume that the discrete-time dynamic of $\lambda_k^i$ is governed by a stochastic process of an unobserved ``state" vector denoted by $\xb_k:=\left[\xb_k^1,\ldots,\xb_k^m\right]$ with $\xb_k^i\in\mathbf{R}^d$:

\begin{equation}\label{eq:State}
    \xb_k = \Psi(\xb_{k-1};\pb)+\eta_k,
\end{equation}
where a function $\Psi$ can be nonlinear, $\pb$ is a fixed parameter vector and $\eta_k\sim N(\mathbf{0},\Qb)$. The conditional intensity, $\lambda_k^i$, is assumed to be a function of the state vector, i.e.,
\begin{equation}\label{eq:observationeq}
    \lambda_k^i = h(\xb_k^i;\qb),
\end{equation}
where the link function of observation, $h$, is usually nonlinear and $\qb$ is a fixed parameter vector. In this setting, the parameter vector is given by the augmented vector $\theta = \left[\pb\enspace\qb\right]$. The so-called complete data likelihood function is the joint probability density $p(\xb_{0:K},\Delta N_{1:K}\mid\theta)$, where $\xb_{0:k}$ denotes the sequence of $\xb_0$ up to $\xb_k$. The maximization problem for ~\eqref{eq:LL} can also be expressed by
\begin{equation}\label{eq:marginalLH}
\widehat{\theta}:=\arg\max_{\theta\in\Theta}~\log\int p\left(\xb_{0:K},\Delta N_{1:K}\mid \theta\right)d\xb_{0:K}.
\end{equation}
\par
We adopt the EM framework to construct an iterative algorithm for the state-space model to avoid a direct integration of the above joint density. The construction of our algorithm follows a similar approach for model identification for the Kalman filter presented in \cite{shumwaystoffer82, Godsill2004MonteCS}.
To this end, we denote the parameter estimate after $\kappa$ iterations by $\theta^{(\kappa)}$. In the EM approach, we have to design a tight lower-bound function (i.e. minorization) that would be more tractable for maximization than the original marginal likelihood function. 
For the current case, a tight lower-bound (or surrogate) function for maximization is given by
\begin{equation}\label{eq:Qfun}
\begin{aligned}
\mathcal{Q}\left(\theta;\theta^{(\kappa)}\right)&=\int p\left(\xb_{0:K} \mid \Delta N_{1:K}, \theta^{(\kappa)}\right) \log p\left(\xb_{0:K}, \Delta N_{1:K} \mid \theta\right) d \xb_{0:K} \\
& =\mathbb{E}\left[\log p\left(\xb_{0:K}, \Delta N_{1:K} \mid \theta\right)\right] .
\end{aligned}
\end{equation}
which represents the E-step of the EM algorithm. The M-step then solves the maximization problem
\begin{equation}
\theta^{(\kappa+1)}:=\arg\max_{\theta\in\Theta}~\mathcal{Q}\left(\theta;\theta^{(\kappa)}\right).
\end{equation}
Under the (first-order) Markovian assumption, we can decompose the surrogate function $\mathcal{Q}\left(\theta, \theta^{(\kappa)}\right)$ by
\begin{equation}\label{eq:Qfunsplit}
\begin{aligned}
& \mathcal{Q}\left(\theta, \theta^{(\kappa)}\right)=Q_0\left(\theta, \theta^{(\kappa)}\right)+Q_x\left(\theta, \theta^{(\kappa)}\right)+Q_{\Delta N}\left(\theta, \theta^{(\kappa)}\right), \\
 &Q_0\left(\theta, \theta^{(\kappa)}\right)=\int p\left(\xb_0 \mid \Delta N_{1:K}, \theta^{(\kappa)}\right)\log p\left(\xb_0 \mid \theta\right)d\xb_0, \\
& =\mathbb{E}\left[\log p\left(\xb_0 \mid \theta\right)\right], \\
& Q_\xb\left(\theta, \theta^{(\kappa)}\right)=\sum_{k=1}^K \int p\left(\xb_k, \xb_{k-1} \mid \Delta N_{1:K}, \theta^{(\kappa)}\right) \log p\left(\xb_k \mid \xb_{k-1},\Delta N_{1:K}, \theta\right) d \xb_k d \xb_{k-1}, \\
& =\sum_{k=1}^K \mathbb{E}\left[\log p\left(\xb_k\mid \xb_{k-1}, \Delta N_{1:K},\theta\right)\right],\\
& Q_{\Delta N}\left(\theta, \theta^{(\kappa)}\right)=\sum_{k=1}^K \int p\left(\xb_k \mid \Delta N_{1:K}, \theta^{(\kappa)}\right) \log p\left(\Delta N_k \mid \xb_k, \theta\right) d \xb_k, \\
& =\sum_{k=1}^K \mathbb{E}\left[\log p\left(\Delta N_k \mid \xb_k, \theta\right)\right] .
\end{aligned}
\end{equation}
\par
To maximize $\mathcal{Q}$, we must assume the availability of $p\left(\xb_0 \mid \theta\right)$, $p\left(\xb_k\mid \xb_{k-1}, \Delta N_{1:K},\theta\right)$, and $p\left(\Delta N_k \mid \xb_k, \theta\right)$.
The initial density $p\left(\xb_0 \mid \theta\right)$ may depend on the model parameter in general, depending on how we would like to generate the initial density for the state. If not, $\mathcal{Q}_0$ can be excluded from the maximization.
The transition density $p\left(\xb_k\mid \xb_{k-1}, \Delta N_{1:K},\theta\right)$ will depend on the model~\eqref{eq:State}. Assuming a normal distribution for $\eta_k$ in~\eqref{eq:State}, $p\left(\xb_k\mid \xb_{k-1}, \Delta N_{1:K},\theta\right)$ is also normal. The likelihood function $p\left(\Delta N_k \mid \xb_k, \theta\right)$ then follows the assumption in~\eqref{eq:PoiLikelihood}.
\par
The expression in \eqref{eq:Qfun} suggests that if we can sample from $p\left(\xb_{0:K} \mid \Delta N_{1:K}, \theta^{(\kappa)}\right)$, we can then estimate all the expectations in~\eqref{eq:Qfunsplit} using the sample paths, which we denote by $\xb_{0:K}^s$. The superscript $s$ stands for ``smoothing" which will be explained below. The efficiency of the EM algorithm in this setting will depend strongly on the design of the path sampling technique. We will use the forward filtering-backward sampling procedure to obtain samples approximately from the joint smoothing distribution \cite{Godsill2004MonteCS}, which is a combination of particle filtering (PF) and backward simulation smoother (BSS) to generate $\xb_{0:K}^s$. 
\par
Particle filtering is a sequential Monte Carlo (SMC) technique for non-linear filtering. In the current application, it can be used to sample $p\left(\xb_{k} \mid \Delta N_{1:k}, \theta^{(\kappa)}\right)$, i.e., only the observation up to the time step $k$ is used to estimate $\xb_{k}$. It enjoys great flexibility but suffers from filtering degeneracy, where most of the sample weights become zero as time increases. A resampling is required to mitigate this issue. We will still, however, employ it in our work for a low-dimensional problem. A brief discussion of PF algorithm is provided in Appendix A. An extensive review of SMC and PF can be found in many review literature, to name a few here, \cite{GordonSalmond93, Doucet2006SequentialMC, Doucet2008ATO,Kitagawa1996MonteCF}.
\par
Suppose that we have obtained the weighted, filtered particle $\left(\xb^{f(i)}_k,\wb^{f(i)}_k\right)$ for $i=1,\ldots,N_f$, where $N_f$ is the number of particles used in PF. The particles approximate $p\left(\xb_{k} \mid \Delta N_{1:k}, \theta^{(\kappa)}\right)$ but the EM algorithm requires $p\left(\xb_{k} \mid \Delta N_{1:K}, \theta^{(\kappa)}\right)$ for any $k=1,\ldots, K$. The BSS uses the filtered particles $\left(\xb^{f(i)}_k,\wb^{f(i)}_k\right)$ to generate the smoothing particles, $\left(\xb^{s(\ell)}_k,\wb^{s(\ell)}_k\right)$ for $\ell=1,\ldots,N_s$, where $N_s$ and $N_f$ can be different. The algorithm for BSS is also provided in Appendix A.  We can then approximate the expectations in~\eqref{eq:Qfunsplit} based on $\left(\xb^{s(\ell)}_k,\wb^{s(\ell)}_k\right)$. Maximization of $\mathcal Q$ in~\eqref{eq:Qfunsplit} to find $\theta^{(\kappa+1)}$ is then carried out numerically. We use the function \verb|fmincon| in \verb|MATLAB| to optimize $\mathcal Q$. A useful by-product of this approach is that the ensemble of sample paths of the conditional intensity $\lambda^i_{1:K}$ can be directly computed from the particles $\left(\xb^{s(\ell)}_k,\wb^{s(\ell)}_k\right)$ at the last iteration of the EM algorithm. The ensemble-based EM method allows for uncertainty quantification of the intensity paths. In the subsequent subsections, we will demonstrate how the ensemble-based EM may be used for some discrete-time Hawkes model that can be represented in a state-space form.

\subsection{Log-Gaussian Cox process (LGCP)}\label{s:LGCP}
We consider a univariate LGCP (i.e. $m=1$) given by
\begin{equation}\label{eq:LGCP1node}
    \begin{aligned}
        x_k&= \underbrace{(1-\omega_1\delta t)x_{k-1}+\omega_1\mu\delta t}_{:=\Psi_x(x_{k-1})}+\epsilon\sqrt{\delta t}\eta_k,\\
        g_k &= \underbrace{(1-\omega_2\delta t)g_{k-1}+\alpha\Delta N_{k-1}}_{:=\Psi_g(g_{k-1})},\\
        \lambda_k &= \exp(x_k)+g_k.
    \end{aligned}
\end{equation}
We assume $\eta_k$ has the standard normal distribution $N(0,1)$. The parameter vector  is $\theta=\left[\mu,\omega_1,\epsilon,\alpha,\omega_2\right]$ and the state variable is $x_k\in\mathbb{R}$. All parameters are assumed to be positive. Hence, ``inhibition" (i.e. $\alpha<0$ ) is not considered in this work. 
\par
We first consider a synthetic experiment with a ground truth $\theta^\ast=\left[1.5, 0.5, 2.5, 0.5, 1.5\right]$ and simulate $\Delta N_{1:K}$ and $\lambda_{1:K}$ for $K=4000$ with $\delta t=0.1$ and initial condition $x_0=1.5$ and $g_0=0$.
We initialize the EM algorithm with parameter vector $\theta^{(0)}=\left[3, 0.25, 1.25, 0.25, 0.75\right]$.
At the $\kappa-$th iteration, the E-step requires a prior sample of the state vector $\xb_0$. We use a prior assumption $x_0\sim N(\mu^{(\kappa)},\epsilon^{(\kappa)}\delta t)$ and set $g_0=0$. The number of particles is $N_f=400$ and we set $N_s=0.25 N_f$ in our experiment. After obtaining the smoothing particle $\xb^{s(\kappa)}$ from E-step (using a combination of PF and BSS as explained in Appendix A), we can evaluate the $\mathcal{Q}-$function in~\eqref{eq:Qfunsplit}. The $\mathcal{Q}-$function (after omitting the terms irrelevant to maximization) has the following form:

\begin{equation}\label{eq:Qeval_LGCP}
    \begin{aligned}
    \mathcal{Q}_0(\theta,\theta^{(\kappa)})&=-\frac{1}{2N_s}\sum_{\ell=1}^{N_s}\frac{\left(x_0^{s(\ell)}-\mu^{(\kappa)}\right)^2}{\delta t}\\
    \mathcal{Q}_\xb(\theta,\theta^{(\kappa)})&=-\frac{1}{2N_s}\sum_{\ell=1}^{N_s}\left[\sum_{k=1}^{K}\frac{\left(x_k^{s(\ell)}-\Psi_x\left(x_{k-1}^{s(\ell)}\right)\right)^2}{\epsilon^2\delta t}+K\log\epsilon\right]\\
    \mathcal{Q}_{\Delta N}(\theta,\theta^{(\kappa)})&=\frac{1}{2N_s}\sum_{\ell=1}^{N_s}\left[\sum_{k=1}^{K}\Delta N_k\log\lambda_k-\exp(\lambda_k)\delta t\right]
    \end{aligned}
\end{equation}
When $K$ is large, the term $\mathcal{Q}_{0}(\theta,\theta^{(\kappa)})$ in~\eqref{eq:Qfunsplit} can be neglected. Furthermore, $\mathcal{Q}_{x}$ depends on only $\mu,\omega_1,\epsilon$ and  $\mathcal{Q}_{\Delta N}$ depends only on $\alpha,\omega_2$. We can then solve the two maximization problems in parallel to obtain $\theta^{(\kappa+1)}$. Maximizing $\mathcal{Q}_{x}$ has a closed-form solution, see Appendix B, and can be readily computed. However, since $\mathcal{Q}_{x}$ and $\mathcal{Q}_{\Delta N}$ are maximized in parallel, numerically maximizing $\mathcal{Q}_{\Delta N}$ becomes a bottleneck to the speed of the algorithm.  For this experiment, we numerically maximize both  $\mathcal{Q}_{x}$ and $\mathcal{Q}_{\Delta N}$ using the interior point algorithm implemented in \verb|fmincon| in \verb|MATLAB|~\cite{fminconByrd}. The constraint optimization is required to ensure the positive values of parameters. 

\par
\begin{figure}[htbp]
    \centering
    \includegraphics[width=\linewidth]{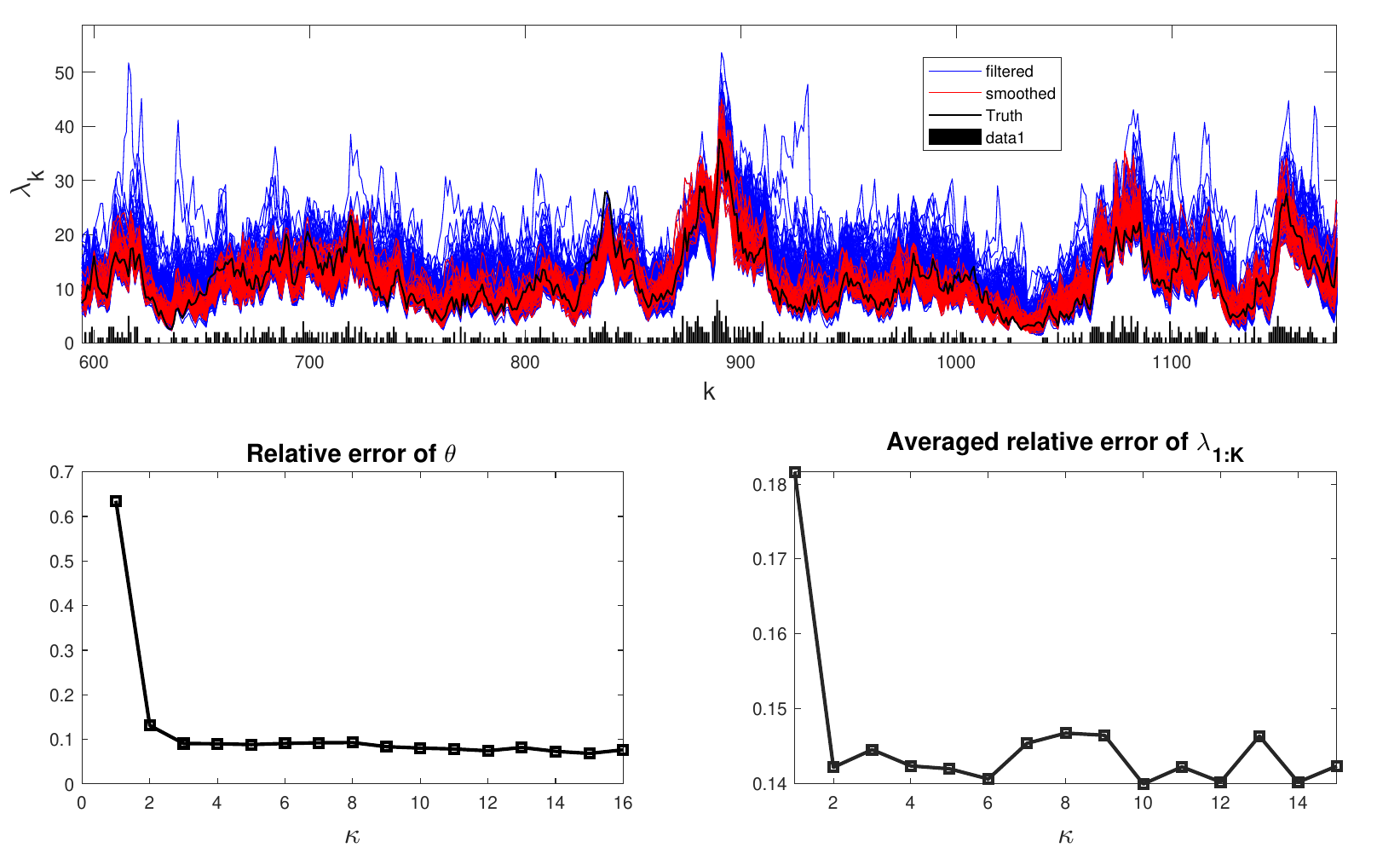}
    \caption{(Top) Comparison between filtered particles $\lambda_{0:k}^{f(\ell)}$, smoothed particles $\lambda_{0:k}^{s(\ell)}$ and ground truth. The simulated data is also shown in the bar plot. (Bottom, left) The progress of the relative error of the parameter vector $\theta$ over 16 EM iterations. (Bottom, right) The propagation of the relative error for the conditional intensity. Note that we have fixed the number of EM iterations to 16; no stopping criteria have been implemented in this example.}
    \label{fig:LGCP1node}
\end{figure}

The results are shown in Figure~\ref{fig:LGCP1node}.

We also compare the filtered intensity $\lambda_{0:k}^{f(\ell)}$ with the smoothed intensity $\lambda_{0:k}^{s(\ell)}$ in Figure~\ref{fig:LGCP1node}, which are computed from $x_{0:k}^{f(\ell)}$, $x_{0:k}^{s(\ell)}$ using the observation equation in \eqref{eq:LGCP1node}, respectively. Figure~\ref{fig:LGCP1node} clearly shows that the smoothed particles have smaller variation (or uncertainty) than that of the filtered particles. The relative error for the parameter vector decays quickly after a few steps and then becomes stable. We compute the relative error for each particle $\lambda_{0:k}^{s(\ell)}$ and report the average relative error in Figure~\ref{fig:LGCP1node}.

\subsection{Logistic LGCP}
We modify the conditional intensity function in \eqref{eq:LGCP1node} to
\begin{equation}\label{eq:LCGP2}
\lambda_k = h\left(\exp(x_k)+g_k)\right),\qquad h(z) = \frac{A}{1+B\exp(-z)}.
\end{equation}
This modification incorporates the upper bound $\lambda_k\leq A$ to the conditional intensity.
The parameter vector  is $\theta=\left[\mu,\omega_1,\epsilon,\alpha,\omega_2,A, B\right]$ and the state vector is $\xb_k\in\mathbb{R}$. The required $\mathcal{Q}-$function is the same as~\eqref{eq:Qeval_LGCP}.
\par
\begin{figure}[htbp]
    \centering
    \includegraphics[scale=0.4]{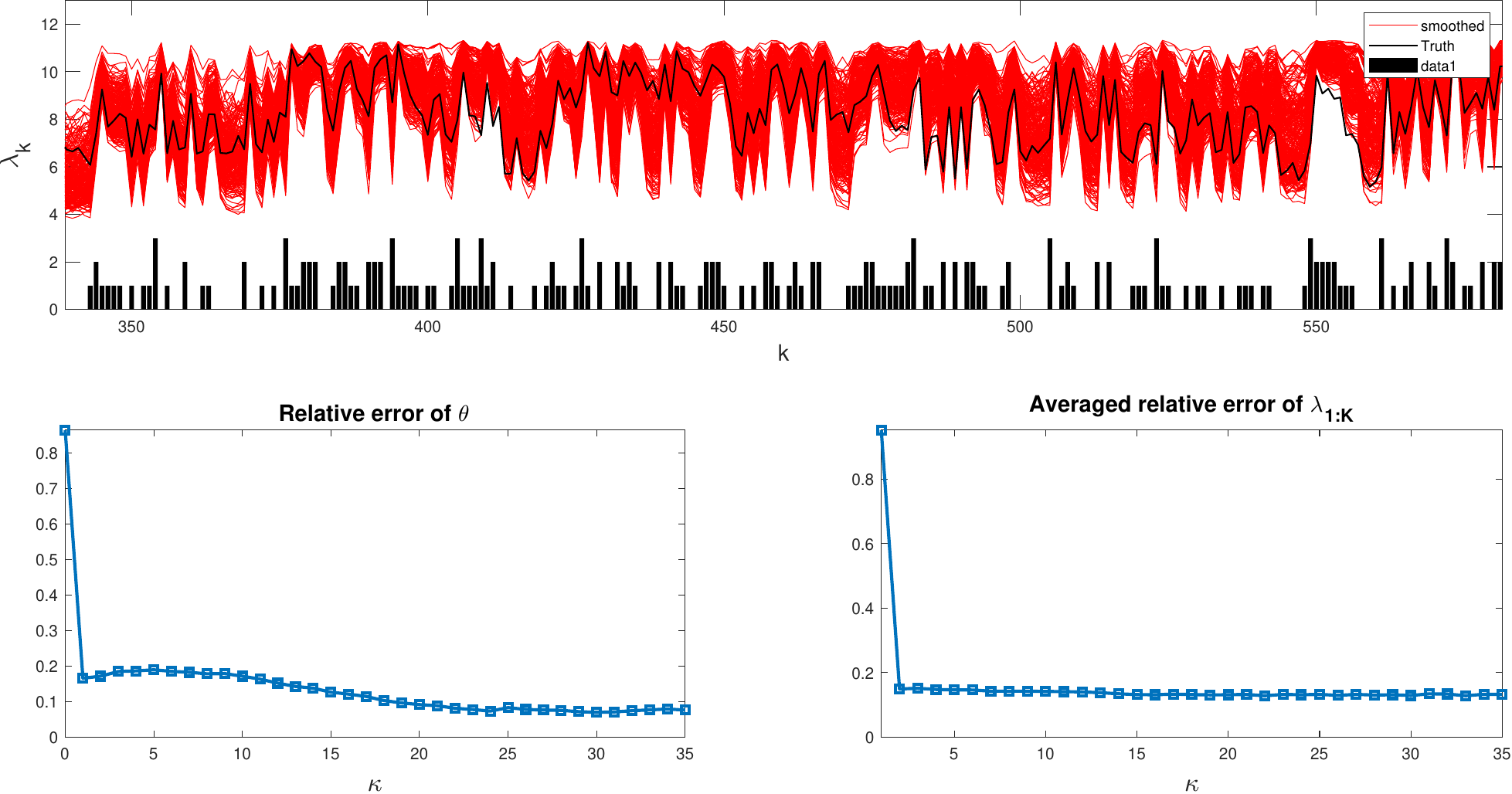}
    \caption{(Top) Comparison between smoothed particles $\lambda_{0:k}^{s(\ell)}$ and ground truth for the Logistic LGCP. The simulated data is also shown in the bar plot. (Bottom, left) The progress of the relative error of the parameter vector $\theta$ over 35 EM iterations where the change in the relative error of $\theta$ is lower than $10^{-5}$. (Bottom, right) The progress of the relative error for the conditional intensity.}
    \label{fig:LGCPlogit}
\end{figure}

We select a ground truth $\theta^\ast=\left[0.5, 0.5, 0.25, 9,  0.5, 12,  4\right]$ and simulate $\Delta_{1:K}$ and $\lambda_{1:K}$ for $K=2000$ with $\delta t=0.1$ and initial condition $x_0=1$ and $g_0=0$. We initialize the EM algorithm with parameter vector $\theta^{(0)}=\left[1, 0.25, 0.5, 4.5, 1, 24, 8\right]$. As shown in Figure \ref{fig:LGCPlogit}, a fast reduction of both relative errors is obtained at the beginning and then slows down as the algorithm converges. This observation is common for MM algorithms in general. In the early stages of the iterations, the parameter estimates are often far from the nearest stationary point (local maximum), resulting in relatively large "gradient components" (in an informal sense). Since the EM algorithm can be roughly interpreted as performing coordinate ascent in the latent or complete-data space~\cite{Langebook16,Mkhardi17,HunterLange04}, these large gradients lead to significant updates in the parameters and, consequently, substantial increases in the likelihood.  

As the algorithm approaches a local maximum (or saddle point), the likelihood surface becomes flatter due to reduced curvature. This causes the parameter updates to shrink, as the re-estimated values in each iteration undergo smaller changes once they are closer to the optimal fit. As a result, the likelihood's improvements diminish with each subsequent iteration.
\subsection{Log-Gaussian Cox process (LGCP) on a small network}
We consider a multivariate extension of LGCP on a network where the links of the network describe an ``influence" structure through a (pairwise) excitation process. We demonstrate the utility of the ensemble EM approach for the multivariate LGCP on the $m$ nodes, denoted by $x^j$ for $j=1,\ldots,m$, that has the following form:
\begin{equation}
    \begin{aligned}
        x^i_{k+1}&=\underbrace{\left[(1-\eta^i)x_k^i+\eta^i\sum_{j\neq i}x_k^j\right](1-\omega_1^i\delta t)+\omega_1^i\mu^i\delta t}_{:=\Psi^i_x(x^i_{k})}+\epsilon^i\sqrt{\delta t}\zeta_k,\\
        g^i_{k+1}&=(1-\omega_2^i\delta t)g_k^i+\sum_{j=1}^m\alpha^{ij}\Delta N_k^j,\quad\lambda^i_{k+1} = \exp(x^i_{k+1})+g^i_{k+1},
    \end{aligned}
\end{equation}
where $\eta^j$ is the diffusion coefficient strength at each node, $\omega_1^j$ and $\omega_2^j$ are the decay rates at each node, $\alpha^{ij}$ are the excitation coupling parameters, $\epsilon^j>0$ and $\zeta_k\sim\mathcal{N}(0,1)$ describes the noise at each node. The $\mathcal{Q}-$function is defined similarly to~\eqref{eq:Qeval_LGCP}.
A (homogeneous) diffusion effect is included in the dynamics of $x^j_k$ and the event-driven excitation process between nodes is incorporated in the dynamic of $g^j_k$. The mutual excitation of $g^j_k$ is driven not only by the count data from the node $j$ itself, but also by all other nodes. Both diffusion and excitation contribute to the increment of the conditional intensity of other nodes in the network in the next time step. 
The state variable  is the vector $\xb_k=\left[x^1_k,\ldots,x^m_k\right]\in\mathbb{R}^m$. 
Note that if there is no diffusion term, we could compute the smoothed path of each $x^j_k$ in parallel for the E-step. 
\par
\begin{figure}[htbp]
    \centering
    \includegraphics[scale=0.6]{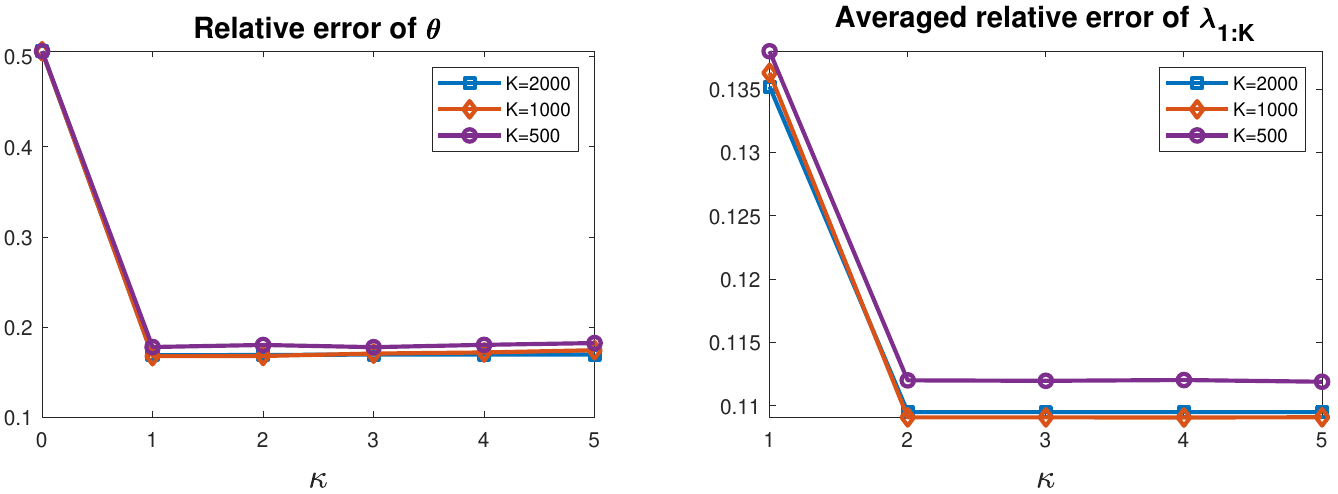}
    \caption{(Left) Relative error of the parameter vector at each EM-step.(Right) Relative error of the conditional intensity at each EM-step. We stop at 5 iterations when the change in the relative error of $\theta$ is lower than $10^{-5}$.}
    \label{fig:LGCP3nodes_err}
\end{figure}
\begin{figure}[htbp]
    \centering
    \includegraphics[width=\linewidth]{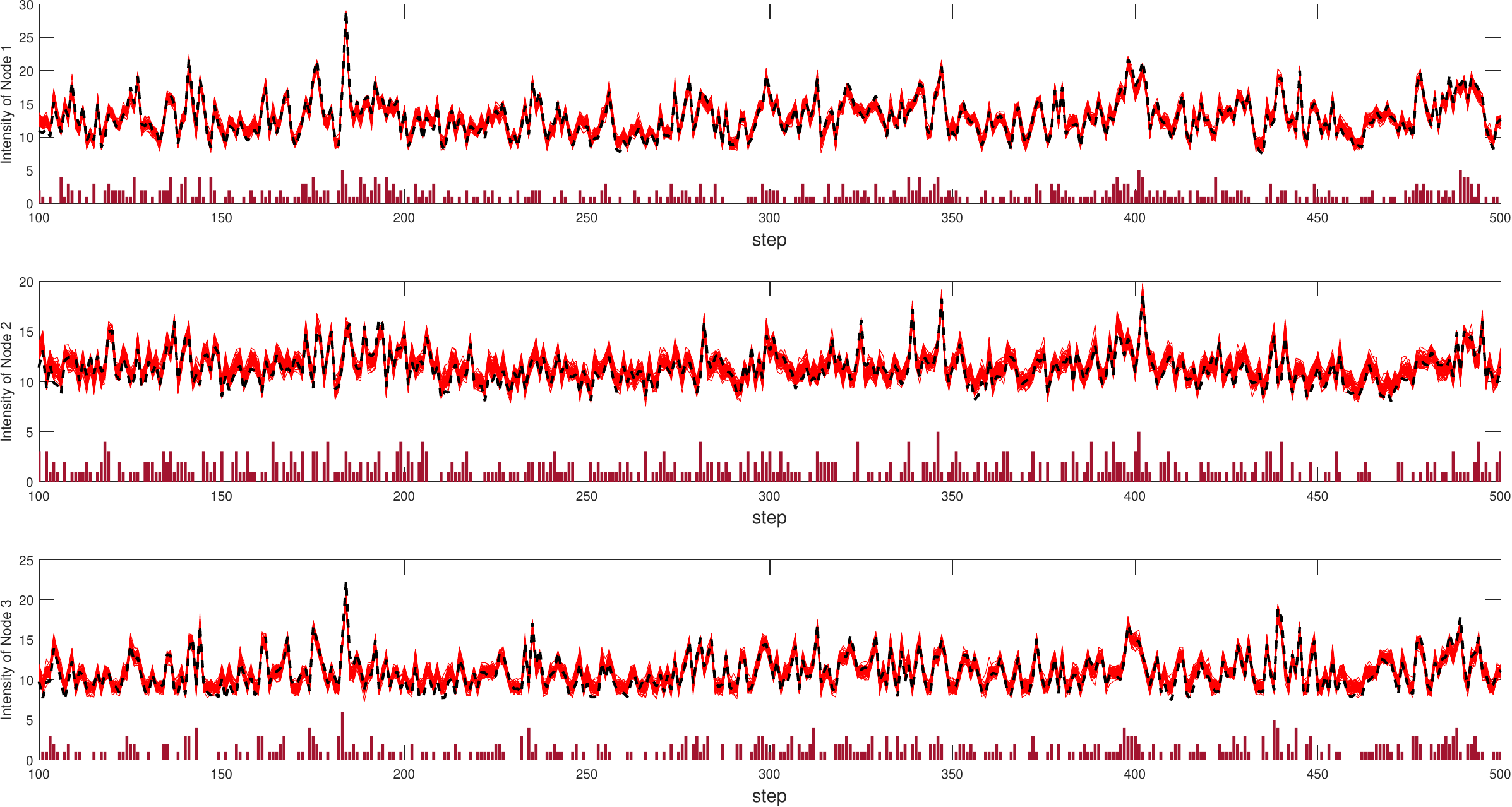}
    \caption{From the top to bottom, the plot shows the smoothed path at the final step of the EM for node $1$ to $3$, respectively. For a clear visualisation, only part of the trajectory is shown at the time step $k=100-500$. The bar plot beneath the intensity shows the simulated count data for each node.}
    \label{fig:LGCP3nodes_smooth}
\end{figure}
\begin{figure}[htbp]
    \centering
    \includegraphics[width=\linewidth]{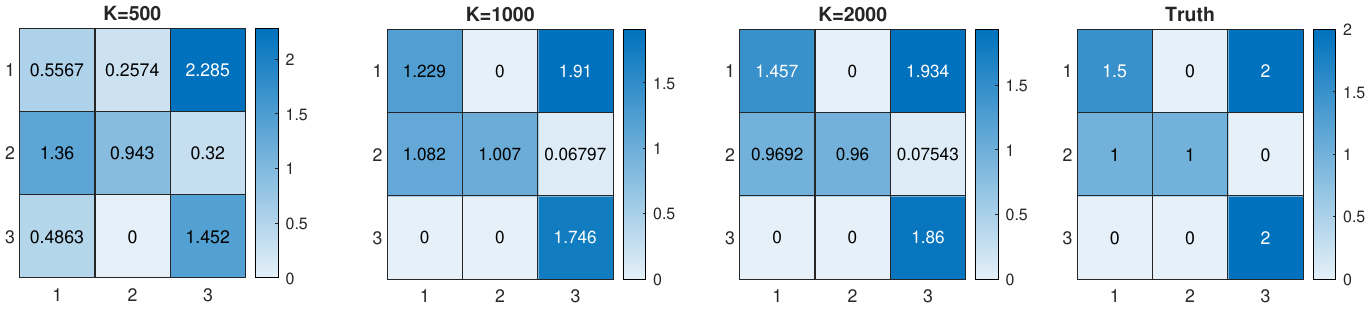}
    \caption{Estimated values of $\alpha^{ij}$ for various data length and the ground truth. The $i-$th row and $j-$ column in the plot indicates $\alpha^{ij}$.}
    \label{fig:LGCP3nodes_NW}
\end{figure}

We consider an experiment with the following parameters: $\delta t=0.1$, $\omega_1^j=0.5/\delta t, \omega_2^j=0.9/\delta t, \eta^j=0.1, \mu^j=0.5, \epsilon^j=0.125$ and $m=3$. The ground truth of the mutual excitation structure $\alpha^{ij}$ is given in Figure \ref{fig:LGCP3nodes_NW}. We test the experiment with different simulated data lengths $K=500,1000,2000$. The initial ensemble for $x^j_0$ is drawn independently from $N(0,5\epsilon^j)$ using $N=600$ particles.  The initial structure of the network is set to $\alpha^{ij}=0.9$, i.e., a fully connected network with a uniform excitation rate of 0.9. The experimental results are shown in Figures \ref{fig:LGCP3nodes_err} to \ref{fig:LGCP3nodes_NW}. The errors of the parameter estimation for various values of $K$ are shown in Figure \ref{fig:LGCP3nodes_err}, all of which exhibit fast error reduction in the first step and then slowly decrease afterwards, similar to the univariate case shown in \S\ref{s:LGCP}. The smoothed path at the final EM step is shown in Figure \ref{fig:LGCP3nodes_smooth} for $K=2000$, which demonstrates a good estimate of the true intensity. The results for the other data lengths are similar. Most importantly, the network structure, which is the main interest of this work, can be accurately captured as shown in Figure \ref{fig:LGCP3nodes_NW} if the data length is sufficiently long enough. Note that we have tested several cases and found similar results when the initial guess of the parameters is ``close enough" to the true parameters in a sense that the stability of the model is sustained; if the initial parameter is not ``close enough" to the true values, the method may fail to converge.

\section{Majorization-Minimization (MM)}\label{sec:MM}
In this section, we focus only on the discrete-time model analogous to the exponential-decay kernel of the multivariate Hawkes process.
In particular, the conditional intensity $\lambda_k^j$ is given by
\begin{equation}\label{eq:dh}
 \lambda^i_{k+1}=\mu^i+(\lambda^i_k-\mu^i)\gamma^i +\sum_{j=1}^m\alpha^{ij}\Delta N^j_k,\quad i=1,\dots,m,
\end{equation}
where $\mu^i>0,\alpha^{ij}\geq0$, $0<\gamma^i<1$ and $\Delta N_{0}^j=0$. The parameter $\mu^i$ represents the baseline rate where the number of events are endogenously generated based on a Poisson distribution with the mean $\mu^i$. We assume that the initial condition is the same as the baseline, i.e., $\lambda^i_{0}=\mu^i$. The parameter $\alpha^{ij}$ for $i\neq j$ models increases in the likelihood to generate more counts for the $i$-th node immediately after observing counts for the $j$-th node. The parameter $\gamma^i$ is the decay rate of $\lambda^i_{k}$ toward the baseline rate. The recursive model~\eqref{eq:dh} can also be rewritten in a closed form by
\begin{equation}\label{eq:integrated}
\lambda^i_{k}=\mu^i+\sum_{\ell=1}^kB_\ell^i\left(\gamma^{i}\right)^{k-\ell}\quad\&\quad B_\ell^i:=\sum_{j=1}^m\alpha^{ij}\Delta N_{\ell-1}^j. 
\end{equation}

This section provides an iterative procedure to minimise the negative log-likelihood function of the discrete-time Hawkes process. In the continuous-time setting, where the timestamp data is available, the branching process can be (artificially) assumed to define the missing data (i.e. immigrant, ancestor or descendent) and the complete likelihood function. However, it is difficult to replicate this idea in the count data setting. Instead, we employ the MM technique to derive an EM-like algorithm for the multivariate Hawkes model driven by count data.
\par
We first present a derivation of the MM algorithm for the case that the decay rate $\gamma^i=\gamma$ is fixed and known. Due to the conditional independence of $\lambda_k^i$ for all $i=1,\ldots,m$ and $k=1,\ldots,K$ given the parameters, we can separately minimise the negative log-likelihood function of each node to estimate $\mu^i$ and $\alpha^{ij}$. For this reason, when estimating parameters for the $i$-th node, we will avoid the notation clutter by omitting the superscript $i$ from the subsequent discussion in this section, e.g., $\alpha^{ij}$ will be written by $\alpha^j$ without confusion. The negative log-likelihood for a given node is 
\begin{equation}\label{eq:NLL}
    \mathcal{L}(\theta):=-\sum_{k=1}^K\log(\lambda_k)\Delta N_k+\delta t\sum_{k=1}^K\lambda_k+\mathcal{C},
\end{equation}
where $\theta=(\mu,\alpha^1\ldots,\alpha^m)$ and $\mathcal{C}$ is a constant. The MM algorithm is an iterative technique that updates the estimation of $\theta^{(n+1)}$ at the $n+1$ iteration by minimising a surrogate function $Q(\theta\mid\theta^{(n)})$. For the minimisation problem, the surrogate function is chosen to be a ``tight" upper bound function so that $Q(\theta\mid\theta^{(n)})\geq\mathcal{L}(\theta)$ for any $\theta$ and $Q(\theta^{(n)}\mid\theta^{(n)})=\mathcal{L}(\theta^{(n)})$. Without loss of generality, we may assume $\delta t=1$ here and ignore the constant $\mathcal{C}$ as well.
\par
By applying Jensen's inequality, a tight upper bound function of $-\log(\lambda_k(\theta))$ in~\eqref{eq:NLL}, denoted by $Q_k(\theta\mid\theta^{(n)})$, can be constructed by
\begin{equation}
    -\log(\lambda_k)\leq Q_k(\theta\mid\theta^{(n)}):=-\frac{\mu^{(n)}}{\lambda_k^{(n)}}\log\left(\frac{\lambda_k^{(n)}}{\mu^{(n)}}\mu\right)-\sum_{l=0}^{k-1}\sum_{j=1}^{m}\frac{\phi_{klj}^{(n)}}{\lambda_k^{(n)}}\log\left(\frac{\lambda_k^{(n)}}{\phi_{klj}^{(n)}}\phi_{klj}\right),
\end{equation}
where $\phi_{klj}:=\alpha^j(\gamma^j)^{k-l-1}\Delta N_k^j$. 
Clearly, we have $Q_k(\theta^{(n)}\mid\theta^{(n)})=-\log(\lambda_k(\theta^{(n)}))$.  Here $\lambda_k^{(n)}$ is a short hand for $\lambda_k^{(\theta(n))}$ and likewise for $\mu^{(n)}$ and $\phi_{klj}^{(n)}$. We now define a tight upper bound function of $\mathcal{L}(\theta)$ by
\begin{equation}
Q(\theta\mid\theta^{(n)})=-\sum_{k=1}^K Q_k(\theta\mid\theta^{(n)})\Delta N_k+\sum_{k=1}^K\lambda_k.
\end{equation}
We obtain the update equations by setting the derivative of $Q$ to zero to yield
\begin{equation}\label{eq:MMupdate1}
      \mu^{(n+1)}=\frac{1}{K}\sum_{k=1}^K\frac{\mu^{(n)}}{\lambda_k^{(n)}}\Delta N_k,\quad\left(\alpha^j\right)^{(n+1)}=\frac{\sum\limits_{k=1}^{K}\sum\limits_{l=0}^{k-1}\frac{\phi_{klj}^{(n)}}{\lambda_k^{(n)}}\Delta N_k}{(1+\gamma)\mathcal{N}^j+\Delta N^j_{K-1}},
\end{equation}
where $\mathcal{N}^j=\Delta N^j_1+\ldots+\Delta N^j_{K-2}$. Note that we make a second-order approximation of small $\gamma^n$ so that $\gamma^n\approx0$ for $n\geq2$ in order to obtain the update equation of $\left(\alpha^j\right)^{(n+1)}$. Notice that the update equations for each parameter are decoupled, so this step can be computed in parallel.

\par
In general, we can also derive MM algorithm that allows the decay rate to be dependent for every node pair. In other words, the parameter vector for the $i$-th node is given by $\theta=\left(\mu^i,\alpha^{i1},\ldots,\alpha^{im},\gamma^{i1},\ldots,\gamma^{im}\right)$. Again, we will omit the superscript $i$ in the algorithm below because of independence of parameters between nodes. Hence, we will write $\theta=\left(\mu,\alpha^{1},\ldots,\alpha^{m},\gamma^{1},\ldots,\gamma^{m}\right)$. To obtain independent update equations for each parameter (similarly to \eqref{eq:MMupdate1}), further work is required to deal with the second term in \eqref{eq:MMupdate1}. Through the Arithmetic-Geometric Mean inequality, the upper bound function  can be obtained by
\begin{equation}\label{eq:MMupper2}
Q(\theta\mid\theta^{(n)})=-\sum_{k=0}^K Q_k(\theta\mid\theta^{(n)})\Delta N_k+K\mu+\sum_{j=1}^m H^j\mathcal{N}^j,
\end{equation}
where 
\begin{equation}
H^j= \frac{\left(\alpha^j\right)^{(n)}}{2\left(1+\left(\gamma^j\right)^{(n)}\right)}(1+\gamma^j)^2+\frac{2\left(1+\left(\gamma^j\right)^{(n)}\right)}{\left(\alpha^j\right)^{(n)}}(\alpha^j)^2.
\end{equation}
Note that we use a second-order approximation of small $\gamma^j$ to obtain the upper bound function \eqref{eq:MMupper2}. By setting the derivative of $Q$ in \eqref{eq:MMupper2} to zero, we obtain the following update equations
\begin{equation}\label{eq:MMupdate_all}
    \begin{aligned}
      \mu^{(n+1)}&=\frac{1}{K}\sum_{k=1}^K\frac{\mu^{(n)}}{\lambda_k^{(n)}}\Delta N_k,\\
    \left(\alpha^j\right)^{(n+1)}&=-B^j+\frac{\sqrt{(B^j)^2+4A^jC^j}}{2A^j},\qquad\left(\gamma^j\right)^{(n+1)}=-D^j+\frac{\sqrt{(D^j)^2+4D^jE^j}}{2D^j},\\
    \end{aligned}
\end{equation}
where 
\begin{equation}\label{eq:MMupdate3}
    \begin{aligned}
      A^j&=\mathcal{N}^j\frac{1+(\gamma^j)^{(n)}}{(\alpha^j)^{(n)}},\qquad B^j=\Delta N^j_{K-1},\qquad C^j=\sum\limits_{k=1}^{N}\sum\limits_{l=0}^{k-1}\frac{\Delta N_k\phi^{(n)}_{klj}}{\lambda_k^{(n)}},\\
      D^j&=\mathcal{N}^j\frac{(\alpha^j)^{(n)}}{1+(\gamma^j)^{(n)}},\qquad E^j=\sum\limits_{k=1}^{K}\sum\limits_{l=0}^{k-1}(k-l-1)\frac{\Delta N_k\phi^{(n)}_{klj}}{\lambda_k^{(n)}}.
    \end{aligned}
\end{equation}
\par
In practice, a regularisation scheme may be required to obtain useful results. For the system with known decay rate, we apply a regularisation only to the baseline rate update in \eqref{eq:MMupdate1} and keep the update equation of the excitation network unchanged. One of the simplest ways to do this is to change \eqref{eq:MMupdate1} to
\begin{equation}\label{eq:reg_knowndecay}
      \mu^{(n+1)}=\frac{1}{K+b}\left(\sum_{k=1}^K\frac{\mu^{(n)}}{\lambda_k^{(n)}}\Delta N_k+a-1\right),
\end{equation}
for some hyperparameters $a,b>0$. This is equivalent to solving the Maximum a posterior (MAP) problem with a gamma prior distribution of $\mu$ with shape parameter $a$ and inverse scale parameter $b$. We will discuss how we choose hyperparameters $a$ and $b$ in the synthetic experiment in the subsequent section.
\par
Similarly, we may regularise the update equations for both $\mu$ and $\gamma^j$ in \eqref{eq:MMupdate_all} by applying a gamma prior distribution for $\mu$ with the hyperparameters $a$ and $b$ and a beta prior distribution for each $\gamma^j$ with hyperparameters $c$ and $d$. Note that for simplicity, we assume the prior distribution of all parameters to be independent and use the same value of the hyperparameters for all $\gamma^j$. The beta distribution is selected to constrain $\gamma^j$ within the desired interval $(0,1)$. The update equation for $\mu^{(n+1)}$ under this regularisation will be the same as \eqref{eq:reg_knowndecay}. However, to update $\gamma^j$, we must solve a quartic polynomial of the following form
\begin{equation}\label{eq:order4}
     -D^jx^4+(D^j+E^j)x^2+(c-d-E^j)x-a=0,
\end{equation}
where $x$ denotes $\left(\gamma^j\right)^{(n+1)}$ and $E^j$ and $D^j$ are defined in \eqref{eq:MMupdate3}. We can either try to solve \eqref{eq:order4} analytically or numerically. In our work, we solve this numerically at every iteration.

\section{ExPKF for the Discrete-time Hawkes model}\label{s:seq_DA}
The MM algorithm in the previous section is designed to estimate the model parameters in batch. Alternatively, we can also develop a sequential procedure to estimate the parameters. This section presents a sequential (second-order) approximation of the posterior density $p(\theta_k\mid\Delta N_{1:k})$. In particular, we are interested in approximating only the mean and covariance matrix associated with $p(\theta_k\mid\Delta N_{1:k})$. Suppose that we have obtained the approximation of the mean and covariance matrix at the time step $k-1$ denoted $\bar{\theta}_{k-1}$ and $\Pcov_{k-1}$, respectively.
Based on this approximation, we assume a prediction model to generate a prior mean, denoted by $\bar{\theta}_{k|k-1}$, and prior covariance, denoted by $\Pcov_{k|k-1}$. Following the derivation in~\cite{Nara18}, a second-order approximation of $p(\theta_k\mid\Delta N_{1:k})$ (called the Extended Poisson-Kalman Filter (ExPKF)) has the following mean $\bar{\theta}_k$ and covariance update $\Pcov_k$ given by
\begin{equation}\label{eq:PKF}
\begin{aligned}
\Pcov_k^{-1} &= \Pcov_{k\mid k-1}^{-1}+\sum_{i=1}^m\biggl[\biggl(\frac{\partial\log\lambda^i_k}{\partial \theta_k}\biggr)\biggl(\frac{\partial\log\lambda^i_k}{\partial \theta_k}\biggr)^{\top}\lambda^i_k\delta t-(\Delta N^i_k-\lambda_k^i\delta t)\frac{\partial^2\log\lambda^i_k}{\partial^2\theta_k}\biggr],\\
\bar{\theta}_k &= \bar{\theta}_{k|k-1}+\Pcov_k\sum_{i=1}^m\biggl[\biggl(\frac{\partial\log\lambda^i_k}{\partial \theta_k}\biggr)(\Delta N^i_k-\lambda^i_k\delta t)\biggr],
\end{aligned}
\end{equation}
where the gradient vector $\frac{\partial\log\lambda^i_k}{\partial \theta_k}$ and Hessian matrix $\frac{\partial^2\log\lambda^i_k}{\partial^2\theta_k}$ are both evaluated at $\bar{\theta}_{k|k-1}$.
\par
The filtering equation~\eqref{eq:PKF} can be used to sequentially approximate the parameters of the model~\eqref{eq:dh}. We will assume that $\gamma^i$ is fixed for a reason that will be explained later; hence there are $m+1$ unknown parameters for each $\lambda^i_k$. To ensure the positivity of the parameters, we will estimate the log-transformed parameter instead,
\begin{equation}
    \theta_k^i:=[\log\mu_k^{i},\log\alpha_k^{i1}\ldots,\log\alpha_k^{im}]^\top.
\end{equation}
Therefore, we have $\theta_k\in\mathbf{R}^{m(m+1)}$. If $\theta_k$ is meant to be a static parameter vector, it is reasonable to assume the following random-walk model,
   $ \theta_k = \theta_{k-1} + \eta_k$,
where $\eta_k\sim N(0,\Qb_k)$. Thus, we have  $\bar{\theta}_{k|k-1}=\bar{\theta}_{k-1}$ and $\Pcov_{k|k-1}=\Pcov_{k-1}+\Qb_k$. Let $S_k^i=\Delta N_k^i+\gamma^i\Delta N_{k-1}^i+\cdots+(\gamma^i)^{k-1}\Delta N_0^i$, which can be recursively computed by $S_{k+1}^i=\gamma^iS_k^i+\Delta N_{k+1}^i$. It can be checked that the gradient vector required by~\eqref{eq:PKF} has the following form
\begin{equation}\label{eq:Dkj}
\frac{\partial\log\lambda^i_k}{\partial \theta_k}=\left[\frac{\partial\log\lambda^i_k}{\partial \theta_k^1},\frac{\partial\log\lambda^i_k}{\partial \theta_k^2},\ldots,\frac{\partial\log\lambda^i_k}{\partial\theta_k^m}\right]^\top,
\end{equation}
where
\begin{equation}
\frac{\partial\log\lambda^i_k}{\partial\theta_k^i}=
\begin{cases}
\frac{1}{\lambda_k^i}\left[e^{\mu_k^i},\enspace S_k^1e^{\alpha_k^{j1}},\ldots,\enspace S_k^m e^{\alpha_k^{jm}}\right], & \text{if } i=j, \\
\boldsymbol{0}, & \text{if } i\neq j.
\end{cases}
\end{equation}
 Thus, only the $i-$th ``block" of $\frac{\partial\log\lambda^i_k}{\partial \theta_k}$ is non-zero. Recall, that $\gamma^i$ is fixed and known. It follows that the Hessian has a simple form, given by
\begin{equation}\label{eq:rank1condition}
\biggl(\frac{\partial\log\lambda^i_k}{\partial \theta_k}\biggr)\biggl(\frac{\partial\log\lambda^i_k}{\partial \theta_k}\biggr)^\top=-\frac{\partial^2\log\lambda^i_k}{\partial^2\theta_k}+\Lambda^{(i)},
\end{equation}
where $\Lambda^{(i)}$ is a diagonal matrix with the diagonal vector $\frac{\partial\log\lambda^i_k}{\partial \theta_k}$. 
Substituting the above results into~\eqref{eq:PKF} yields
\begin{equation}\label{eq:invcovrank1}
\Pcov_k^{-1} = \Pcov_{k\mid k-1}^{-1}+\sum_{i=1}^m\Delta N^i_k\left(\frac{\partial\log\lambda^i_k}{\partial \theta_k}\right)\left(\frac{\partial\log\lambda^i_k}{\partial \theta_k}\right)^\top+\left(\lambda_k^i\delta t-\Delta N^i_k\right)\Lambda^{(j)}.
\end{equation}
With the form in \eqref{eq:invcovrank1}, a rank-1 update can be efficiently used to compute $\Pcov_k^{-1}$. Also, if $\Pcov_{k\mid k-1}$ has a block-diagonal form where each block corresponds to the parameters of each node, $\Pcov_k^{-1}$ will also have the same block-diagonal structure where the $i-$th block corresponds to the parameters of the $i$-th node. Therefore, by ensuring that $\Pcov_{k\mid k-1}$ has the same block-diagonal structure, $\Pcov_{k}$ will inherit the same block-diagonal structure. Consequently, the update system~\eqref{eq:PKF} can be implemented for each node in parallel. To this end, we will always enforce the block-diagonal structure to $\Pcov_0$ and $\Qb_k$ in all of our numerical experiments.

\section{Synthetic data tests}\label{s:results}
\subsection{Test Experiment}\label{sec:test exp1}
We set up this experiment to generate synthetic test data for three different scenarios based on  \eqref{eq:dh}. The true network has $m=9$ nodes with the following baseline rates: $\mu^1=5, \mu^2=4.6, \mu^3=4.2, \mu^4=0.5, \mu^5=0.46, \mu^6= 0.42, \mu^7=0.38, \mu^8=0.34,$ and $\mu^9=0.3$.

\begin{figure}[htb]
    \centering
    \includegraphics[width=\linewidth]{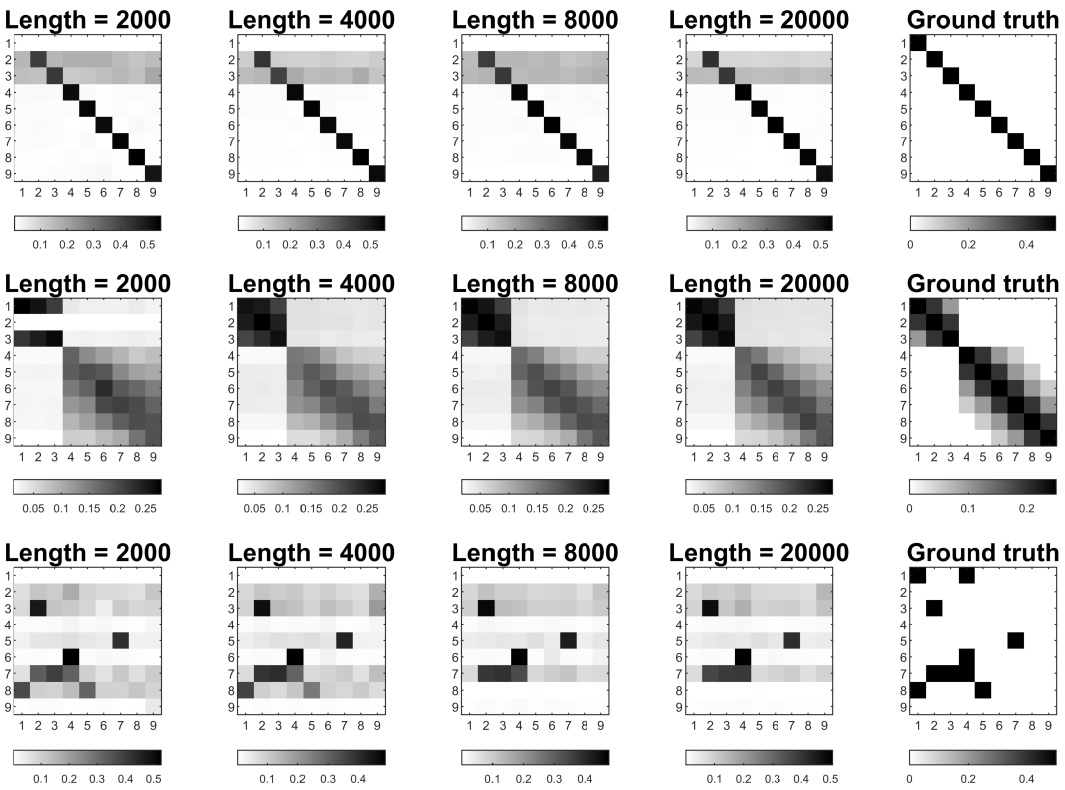}
    \caption{The estimated network structure obtained from the MM algorithm \textbf{without} regularisation for three different ground truths. The length of the data is varied to demonstrate the convergence to the ground truth.}
    \label{fig:9nodesMMnoreg}
\end{figure}

We assume $\gamma^i:=\gamma=0.175$ for all nodes. Note that the model \eqref{eq:dh} is stable if the magnitude of the largest eigenvalue of $\delta t(1-\gamma)^{-1}\Ab$ is less than 1, where $\Ab$ is a matrix with $\alpha^{ij}$ entries on the $i-$th row and $j-$th column. The value of $\gamma=0.175$ is selected so that the model \eqref{eq:dh} is stable for all three different ground truths of the excitation matrix, $\Ab$, examined in this experiment. The structures of three different ground truths are shown in Figure~\ref{fig:9nodesMMnoreg} representing the different scenarios: only self-excitation (top row), localised excitation (middle row), and random excitation structure (bottom row). We generate the test data with $\delta t=1$ for various data lengths, $K=2000, 4000, 8000, 20000$. A test data is simulated by running the model \eqref{eq:dh} and sampling $\Delta N_k^i$ from a Poisson distribution with the mean rate $\lambda_k^i\delta t$ with the initial condition $\lambda_0^i=\mu^i$. 
\par
\begin{figure}[htb]
    \centering
    \includegraphics[width=\linewidth]{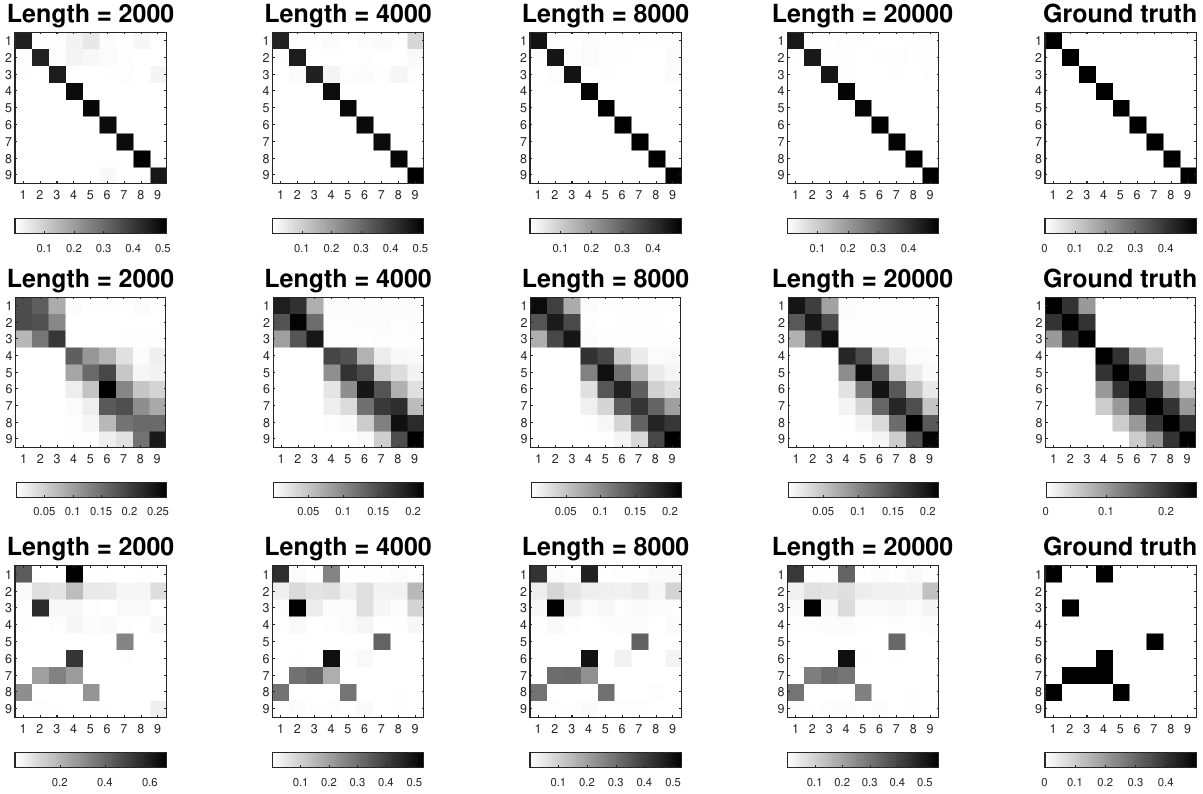}
    \caption{The estimated network structure obtained from the MM algorithm \textbf{with} regularisation for three different ground truths. The length of the data is varied to demonstrate the convergence to the ground truth.}
    \label{fig:9nodesMMregbeta}
\end{figure}

The results for MM method with and without regularisation are shown in Figures \ref{fig:9nodesMMregbeta} and \ref{fig:9nodesMMnoreg}, respectively. For any $i-$th node, the hyperparameter for the regularised MM algorithm is set to $a = 0.5\overline{N^i}K$, where $\overline{N^i}$ is the average count of the $i-$th node over $K$ time steps and $b=K$. This is equivalent to choosing the gamma prior distribution with mean $0.5\overline{N^i}$ (half of the total count for the given node) and variance $0.5\overline{N^i}/K$. Although the selection of these prior parameter values may be arbitrary in general, we based our choice of the prior mean on the reasoning that the baseline rate should be lower than the data average due to the creation of certain events by excitation. In addition, the variance is sufficiently small to ensure that the prior information, or regularisation, is not dominated by the sample size. We also set $c=2.5K$ and $d=10.25K$ for the beta prior distribution, which gives the mean $1/6$ and variance in the order of $o(1/K)$. We make this selection to prolong the impact of excitation by avoiding a decay rate that is too close to 1. This selection is again arbitrary. For the remainder of our work, we will utilize this approach to select prior information for the MM algorithm. Our results, as demonstrated in Figures \ref{fig:9nodesMMregbeta} and  \ref{fig:9nodesMMnoreg}, clearly illustrate that the use of regularization enables the algorithm to produce significantly more accurate outcomes that closely approximate the ground truth.

\par
\begin{figure}[htb]
    \centering
    \includegraphics[width=\linewidth]{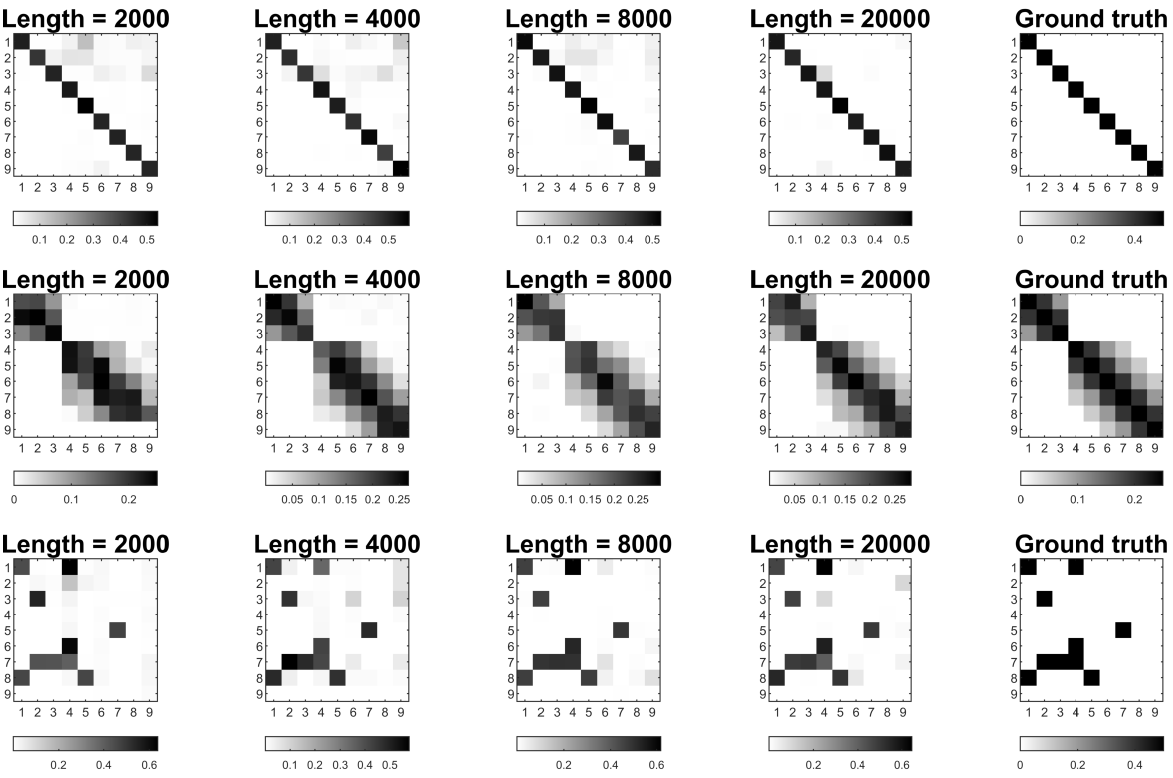}
    \caption{The estimated network structure obtained from the ExPKF algorithm using a true decay rate. The length of the data is varied to demonstrate the convergence to the ground truth.}
    \label{fig:9nodesExPKF}
\end{figure}
\begin{figure}[htb]
    \centering
    \includegraphics[scale=0.5]{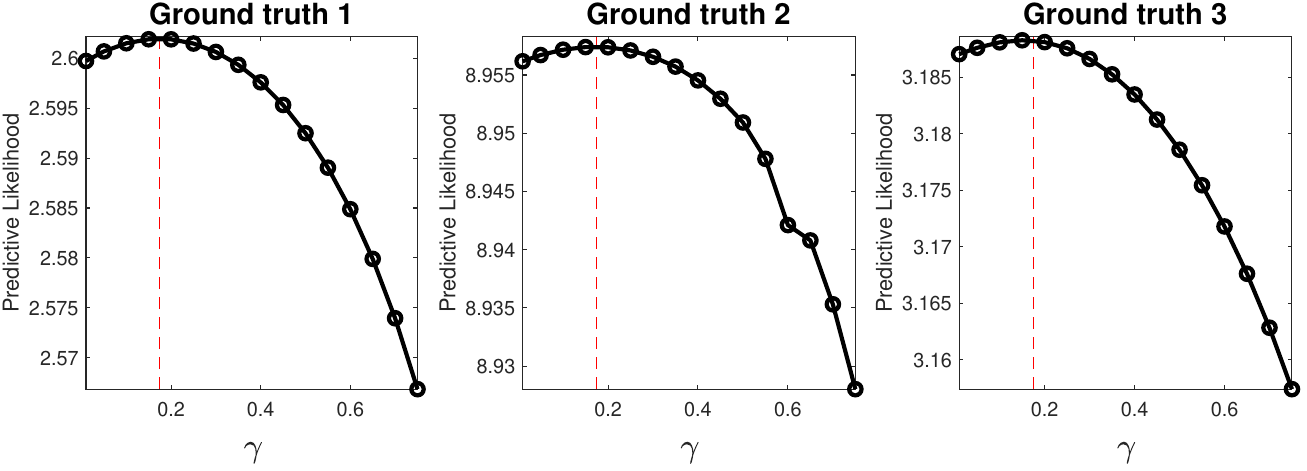}
    \caption{Comparing the predictive log likelihood for various values of decay rate, $\gamma$. The vertical line indicates the true value of the decay rate. The data with length of 20000 is used to produce this result.}
    \label{fig:9nodesLL}
\end{figure}

We use the same simulated data to test ExPKF. Recall that we must assume a fixed decay rate for ExPKF to achieve efficient algorithms via the rank-1 update. We will discuss this issue later in this section. We set $\Pcov_0=10^{-4}\Ib$ and $\Qb_k=10^{-5}\Ib$ for all $k=1,\ldots,K$. The initial guess of the $\mu^j$ is set to be half of the data average of the $j-$th node. Figure \ref{fig:9nodesExPKF} illustrates that the ExPKF algorithm can achieve accuracy comparable to that of the MM algorithm with regularization in capturing network structures. However, ExPKF requires a known decay rate value for all nodes, which was used in this experiment. In practice, decay rates may vary and be unknown for different nodes. To overcome this issue, we explore a method to identify the optimal decay rate $\gamma^i$, assuming a uniform decay rate for all nodes. Specifically, we perform a one-dimensional maximization based on the average log predictive probability \eqref{eq:PoiLikelihood}. To calculate the average log predictive probability, we use the parameter estimate obtained from the previous time step ($\bar{\theta}_{k-1}$) and evaluate \eqref{eq:PoiLikelihood} at time step $k$, averaging over all time steps. Figure \ref{fig:9nodesLL} demonstrates that maximizing the predictive probability yields the optimal decay rate.
\par
\begin{figure}[htb]
    \centering
    \includegraphics[width=\linewidth]{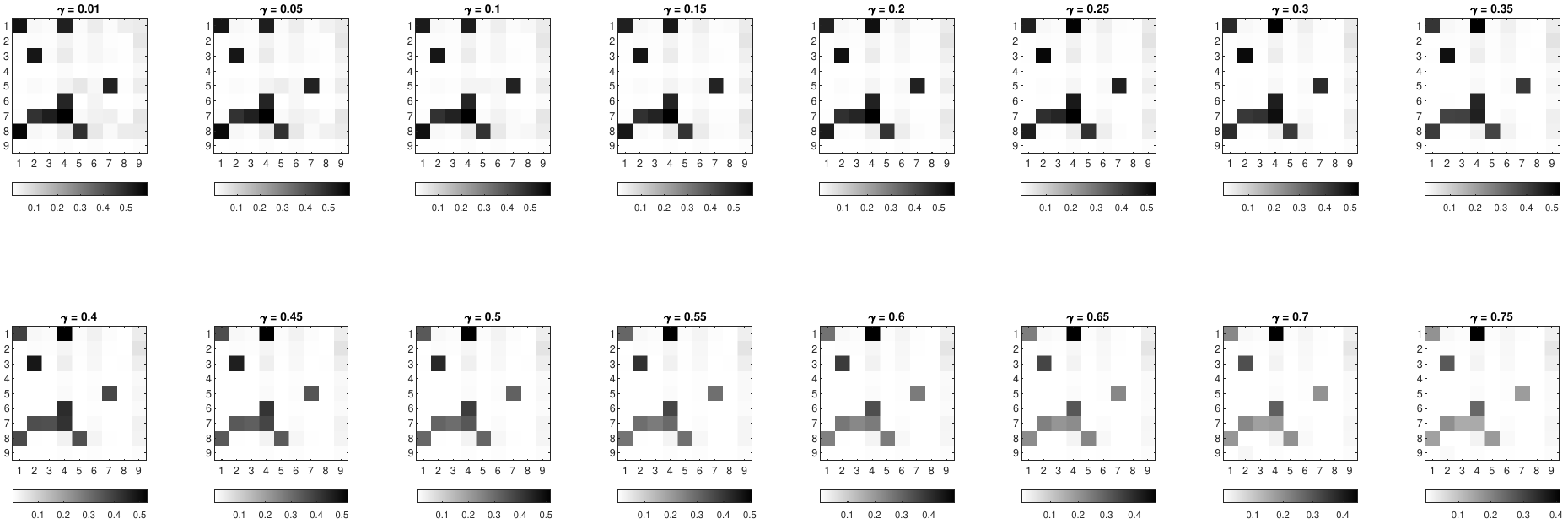}
    \caption{Comparing the network structure for the ground truth 3 when the value of decay rate, $\gamma$, is incorrectly specified.}
    \label{fig:9nodeswrongdecay}
\end{figure}

Moreover, the network structure appears to be highly robust to parameter misspecification, as shown in Figure \ref{fig:9nodeswrongdecay} for the ground truth 3 scenario. Despite the presence of spurious links caused by incorrect decay rates, the primary network structure closely resembles the true structure. Similar results are observed for the other ground truths, although they are not presented in this study.

\par

\subsection{Estimation under model misspecification}
In this experiment, we will generate test data from an agent-based model (ABM) on a set of nodes featuring an excitation network structure, which will then be estimated within the Hawkes model; hence, a model misspecification problem. We adopt a model inspired by the ABM in~\cite{short2008} to simulate the random movement of an ``agent" between nodes through the network edges. During a time interval $(t, t + \delta t)$, an agent located at a node $s$ creates a number of ``events" independently, following a Poisson probability with mean $A_s(t)\delta t$. The total number of events generated at $s$ during $(t, t + \delta t)$ is denoted by $E_s(t)$. The discrete-time dynamic of $A_s(t)$ is given by $A_s(t) = A_s^0 + B_s(t)$, where $A_s^0$ is a static, node-dependent baseline rate, and $B_s(t)$ is a dynamic component that follows the rule:
\begin{equation}
    B_s(t+\delta t)=\left[(1-\eta_s)B_s(t)+\eta_s\sum_{s'\sim s}B_{s'}(t)\right](1-\omega_s\delta t)+\sum_{s'\sim s}w(s,s')E_{s'}(t).
\end{equation}
The interpretation of these parameters is listed below:
\begin{itemize}
    \item $0<\omega_s<1$ is the node-dependent decay rate;
    \item $0<\eta_s<1$ is the  node-dependent diffusion rate;
    \item $w(s,s')\geq0$ defines the strength of the event-driven excitation rate that the node $s'$ has on the node $s$ and we write $s'\sim s$ if $w(s,s')>0$.
\end{itemize}
If an agent generates one or more events, the agent will be removed from the simulation. Otherwise, the agent will move from a node $s$ to $s''$ such that $w(s,s'')>0$ based on a discrete probability distribution
\begin{equation}
    q(s,s''; t):=\frac{A_{s''}(t)}{\displaystyle\sum_{s'\sim s}A_{s'}(t)}.
\end{equation}
Prior to a simulation of the next time step, new agents will be independently created for each node according to a Poisson distribution with mean $\Gamma_s\delta t$ where $\Gamma_s$ is a fixed parameter. The key difference between the Hawkes and ABM models is that the Hawkes model's network structure provides only data-driven excitation, whereas the ABM network structure determines the agents' excitation, diffusion, and probabilistic movement of agents.
\par
The focus of this experiment is on the structure of $w(s,s')$. To this end, we define an ``influence" matrix $\boldsymbol{W}$ such that its $s-$th row and $s'-$th column entry is $w(s,s')$, and reconstruct the pattern of non-zero elements of $\boldsymbol{W}$ using the ExPKF and MM methods. It is important to note that no ground truth is available for this experiment. Although the dynamics of the ABM exhibit similarities to the Hawkes process in terms of influence effects, diffusion effects are also present in the ABM but absent from the Hawkes process. Nevertheless, we anticipate that the influence structure of the Hawkes model will be akin to that of the ABM when fitting the Hawkes model with ABM-simulated count data.
\par
\begin{figure}[htbp]
    \centering
    \includegraphics[width=\linewidth]{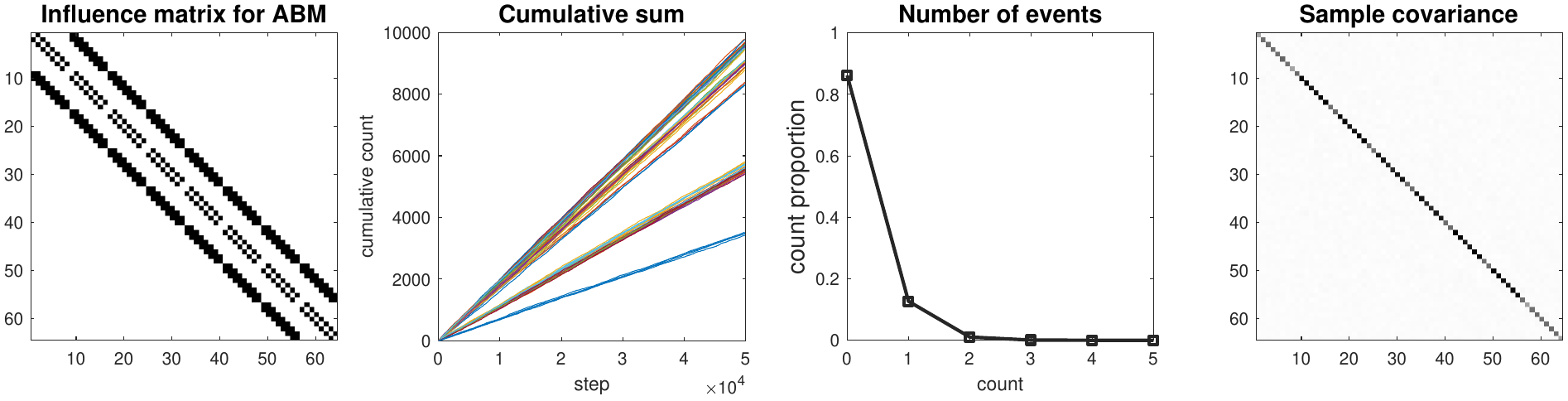}
    \caption{Data simulated from ABM model. The left most plot illustrates the influence matrix $\boldsymbol{W}$. The second plot from the left shows the cumulative number of events for all nodes. The third plot is the frequency distribution of the count. The right most plot is the sample covariance matrix.}
    \label{fig:dN64ABM}
\end{figure}

We generated data of various lengths based on the same influence matrix pattern, represented by $\boldsymbol{W}$ with 64 nodes, as illustrated in Figure~\ref{fig:dN64ABM}. The influence structure is sparse and irreducible, and we considered two cases. In Case 1, we set all non-zero influences to $w(s,s')=3$ for $s,s'\in{1,\ldots,64}$, with $B_s(0)=0$ for all $s$ and $\delta t=0.05$. Other static parameters were spatially uniform: $\omega_s=5$, $\eta_s=0.25$, and $\Gamma_s=3$, for all $s$. Additional information on the simulated data for Case 1 is provided in Figure~\ref{fig:dN64ABM}. The cumulative counts of the data can be grouped into three distinct categories based on the unique values of the row sums of $\boldsymbol{W}$. Most time intervals exhibit either no events or a single event, and the sample covariance of the count time-series indicates little correlation among the nodes.

\par
\begin{figure}[htb]
    \centering
    \includegraphics[scale=0.5]{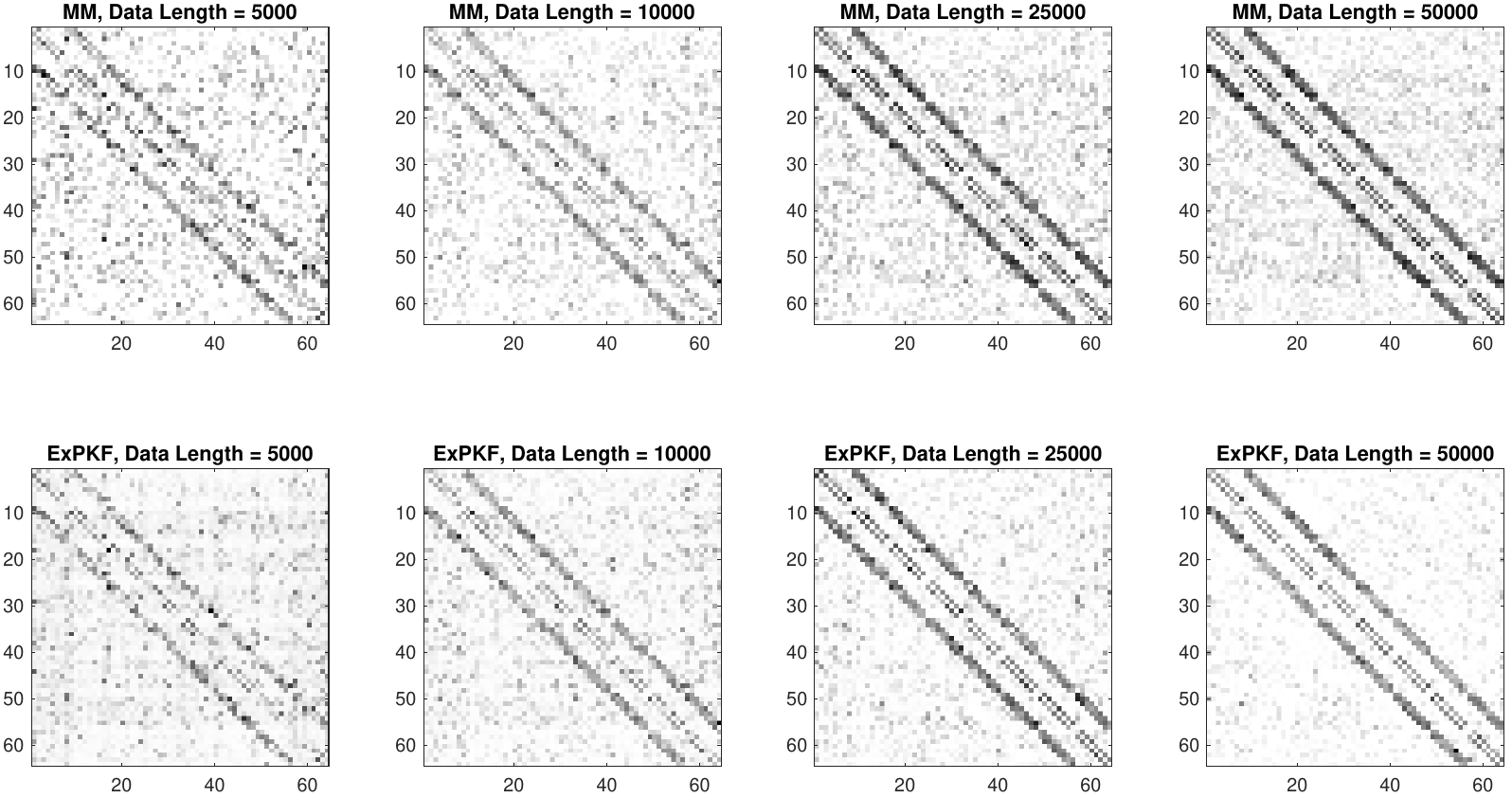}
    \caption{MM and ExPKF estimation of the influence matrix using different data length for Case 1, which has a strong excitation effects.}
    \label{fig:ABMmatrix}
\end{figure}

We use the simulated data to test ExPKF and MM. For ExPKF, we set $\Qb_k=10^{-5}\Ib$ for all $k$ and $\mathbf{P}_0=10^{-4}\Ib$ for all nodes. Again, we set the initial value $\bar{\theta_0}$ in the same manner as done with the previous experiment in Section~\ref{sec:test exp1}. Figure~\ref{fig:ABMmatrix} shows the estimated $\boldsymbol{W}$ based on ExPKF and MM. Interestingly, both methods can correctly reconstruct the pattern of the influence matrix $\boldsymbol{W}$ despite the model misspecification. The results improve with longer data series. Notably, ExPKF produces a result with less ``noise" in the part that is supposed to have zero influence. 
\par
\begin{figure}[htb]
    \centering
    \includegraphics[scale=0.4]{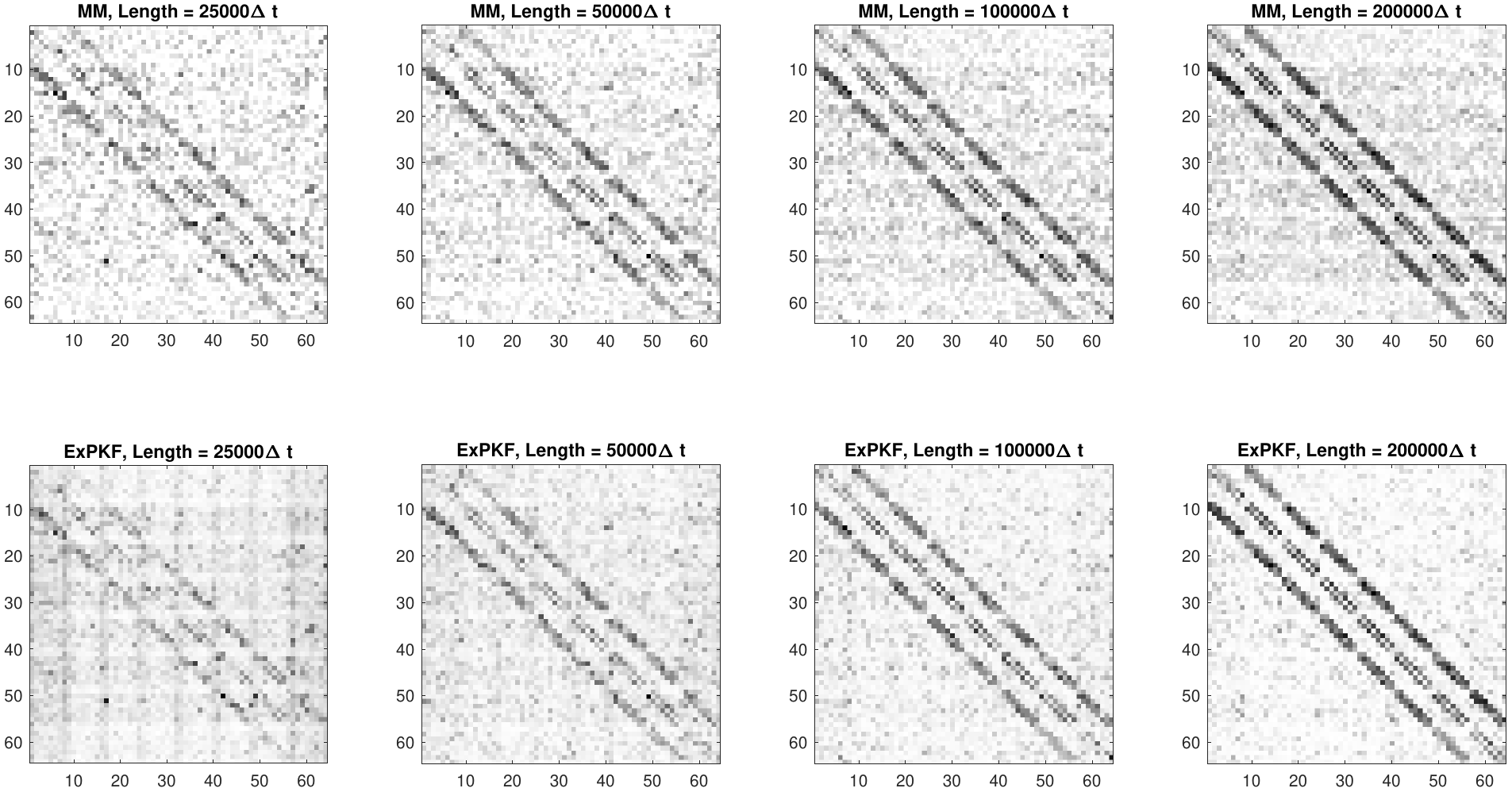}
    \caption{MM and ExPKF estimation of the influence matrix using different data length for Case 2, which has a weak excitation effect but strong diffusion.}
    \label{fig:ABMmatrix_dataNew}
\end{figure}

We also examine Case 2 where we change the non-zero influences to $w(s,s')=0.5, \eta_s=1, \Gamma_s=0.5$, keeping all the other parameter values the same. While the network pattern in Case 1 is manifested mostly through the excitation process (i.e. larger values of $w(s,s')$ and smaller values of $\eta_s$ and $\Gamma_s$), the Case 2 has a weak excitation and generation rate of the new agents but increased diffusion. This change would make it more difficult to detect the influence pattern. Nonetheless, the network structure can still be detected as displayed in Figure~\ref{fig:ABMmatrix_dataNew}. However, the data length required to achieve a good result has to be longer than that of the Case 1; note that the average number of counts per time step is 0.15 for Case 1 but only 0.025 for Case 2.
\par
\begin{figure}[htbp]
    \centering
    \includegraphics[scale=0.425]{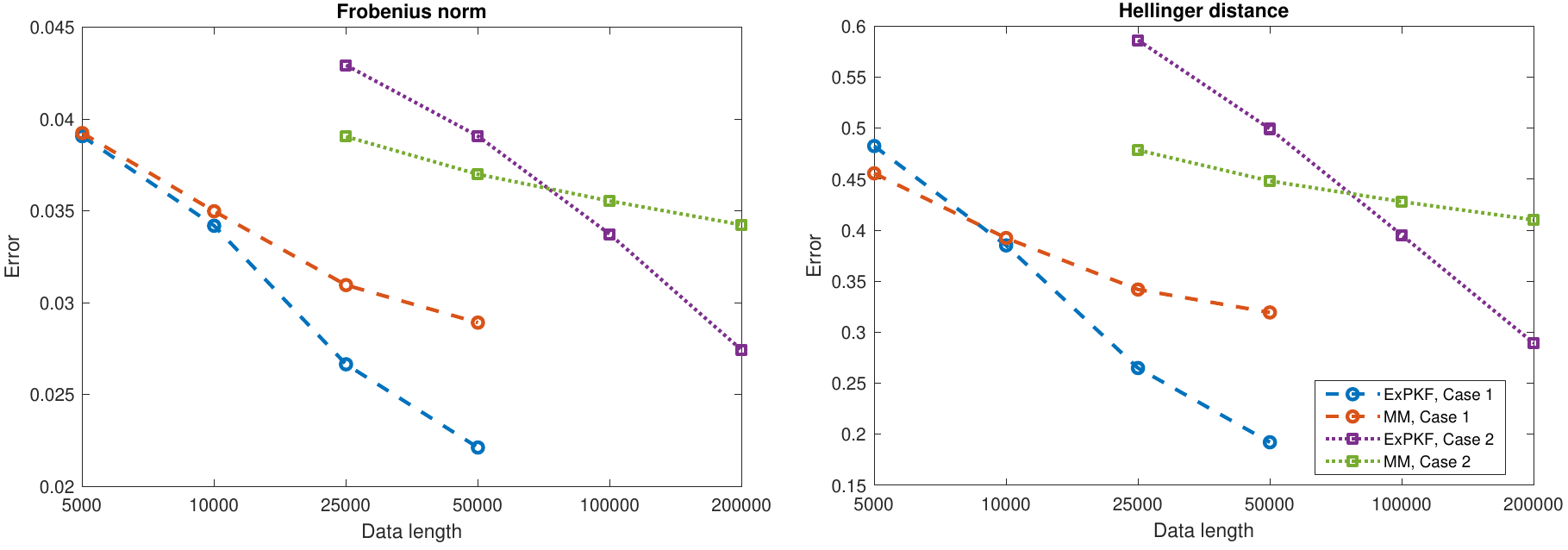}
    \caption{MM and ExPKF Error for different data lengths. (Left) Frobenius norm. (Right) Hellinger distance. Two cases are presented: Case 1 (plotted with the circle markers) corresponds to Figure~\ref{fig:ABMmatrix} and Case 2 (plotted with the square markers) corresponds to Figure~\ref{fig:ABMmatrix_dataNew}.}
    \label{fig:ABMerror}
\end{figure}

The errors based on the Frobenius norm and the Hellinger distance for different data lengths are analysed in Figure~\ref{fig:ABMerror}. When computing both error measures, we normalise $\boldsymbol{W}$ so that the sum of all elements is 1. This is necessary since we have no numerical ground truth to compare against and should evaluate the error based only on the network structure. Both error measures suggest a small improvement of the ExPKF results, which is consistent with the visualisation of Figure~\ref{fig:ABMmatrix}.

\section{Email network data}\label{s:email}
\subsection{Small email network}
\begin{figure}[htb]
    \centering
    \includegraphics[scale=0.5]{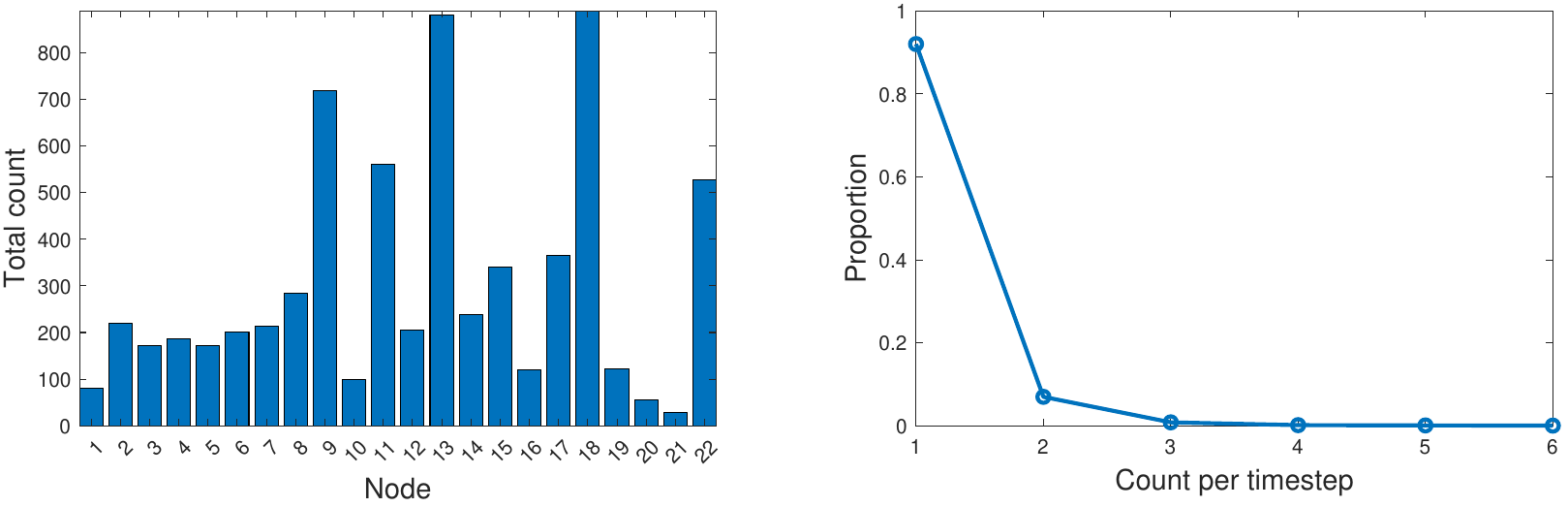}
    \caption{The total count of emails sent by each node (Left) and proportion of non-zero counts (Right).} 
    \label{fig:Ikedatainfo}
\end{figure}

In this section, we analyse the Ikenet dataset consisting of log files from email transactions between 22 anonymized officers at West Point Military Academy over a one-year period, which is available to download via \url{https://github.com/naratips/Ikenet.git}. The dataset contains the time-stamps of outgoing emails and their corresponding receivers. Table~\ref{tab:Ike} displays the top 9 sender-receiver pairs in the dataset, ranked by the number of out-going emails. The data clearly highlights the overwhelmingly large amount of mutual email correspondence for the pairs $(9,18)$ and $(11,22)$. Previous studies on this dataset have utilized information about both the sender and recipients of emails \cite{FOX16,zipkin_schoenberg_coronges_bertozzi_2016}. The Hawkes model with the exponential decay rate was used in~\cite{FOX16} where the rate of sending out an email for a given node is driven by the events of emails received by the given node. For our experiment, we will focus solely on information about the outgoing emails. Therefore, the ``influence" in our analysis can be interpreted as the effect of the number of emails sent out by other nodes on the rate of sending out emails (without any knowledge of the recipients). To demonstrate the methods in terms of count data, we aggregated the timestamp data of outgoing emails into a time-series of count data with a uniform temporal interval of $dt=0.1$days. Figure~\ref{fig:Ikedatainfo} shows the total number of counts for each node and the proportion of non-zero counts per time step.
Note that despite having the number of outgoing emails as large as node 18, node 13 is not among the top pairs $(9,18)$ and $(11,22)$.
\begin{table}[ht!]
\centering
\begin{tabular}{|c|c|c|c|c|c|c|c|c|c|}
\hline
sender&18 & 9 & 22 & 11 & 15 & 8 & 18 & 13 & 18 \\
\hline
receiver&9 & 18 & 11 & 22 & 13 & 18 & 8 & 17 & 22\\
\hline
$\%$ of total &6.95 & 5.97 & 3.97 & 3.01 & 1.96 & 1.89 & 1.87 & 1.78 & 1.75 \\
\hline
\end{tabular}
\caption{Top 9 sender-receiver for the Ikenet data ranked by the total number of outgoing emails.}
\label{tab:Ike}
\end{table}

\par
\begin{figure}[htb]
    \centering
    \includegraphics[scale=0.75]{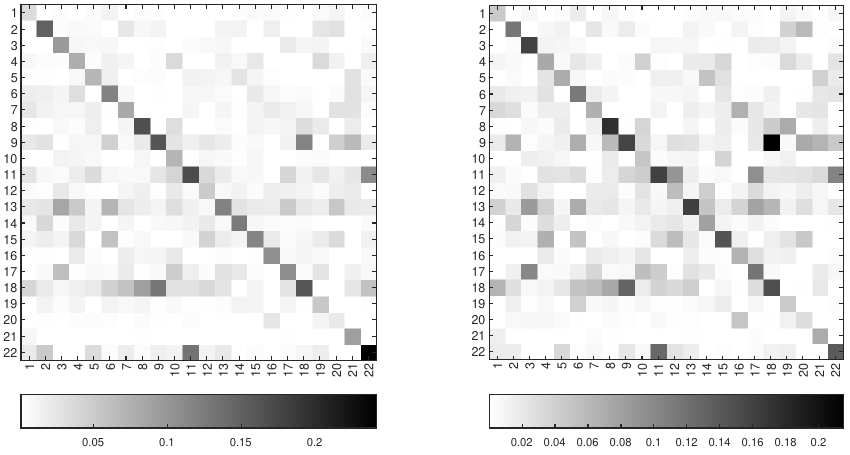}
    \caption{Influence networks associated with the 22-node Ikenet email data constructed by MM algorithm (Left) and ExPKF algorithm (Right).} 
    \label{fig:Ikenw}
\end{figure}

As shown in Figure~\ref{fig:Ikenw}, networks constructed by MM and ExPKF are very similar. By comparing the dominant connections in the network with Table~\ref{tab:Ike}, we can see that the influence network highlights the top sender-receiver pairs $(9,18)$ and $(11,22)$ even though no knowledge of the email recipient network is used in the experiment. 

\subsection{Large email network}
\begin{figure}[htb]
    \centering
    \includegraphics[width=\linewidth]{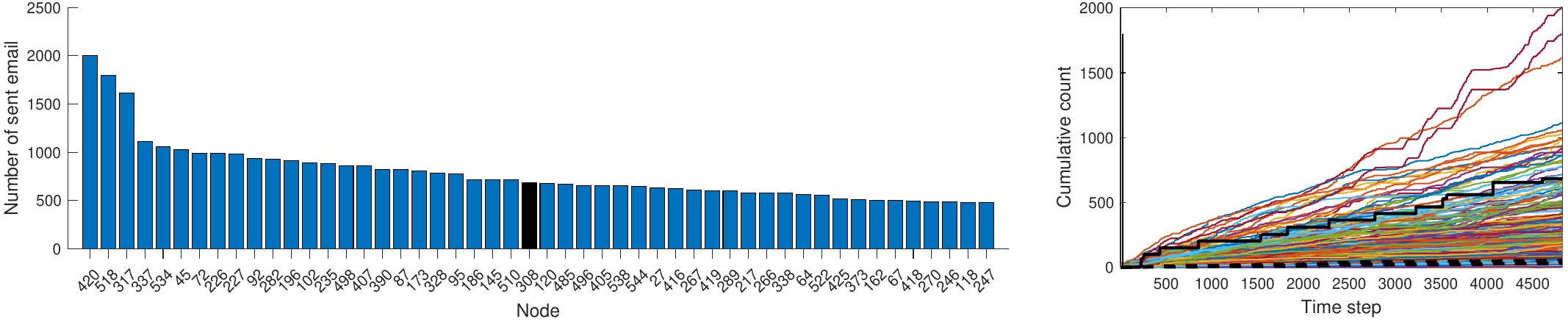}
    \caption{(Left) Histogram of the top 50 users by number of emails sent. The black bar is associated with node 308, which is identified by ExPKF as the most influential node. (Right) The cumulative counts of all nodes. The node 308 has the cumulative counts shown in a black solid curve. The cumulative counts of the nodes influenced by node 308 are plotted in the black dash curves.}
    \label{fig:Email_countdata}
\end{figure}

In this section, we carry out an experiment on a real-world (anonymised) email timestamp data similar to the previous section but at a much larger size. The original data can be found from the following link: \url{https://snap.stanford.edu/data/email-Eu-core-temporal.html}. However, we focus only on the outgoing emails and we ``cleaned up" the data by removing a continuous period of extremely low count due to missing data, weekends and holidays. The cleaned-up data has 545 ``nodes" and 61821 intervals (each interval is one hour long) in total with approximately $97\%$ of zero counts, $2\%$ of one count per interval and the rest of the data has more than one count. Figure~\ref{fig:Email_countdata} shows the top 50 nodes with the highest number of emails sent and the cumulative count for all nodes. 
\par
\begin{figure}[htb]
    \centering
    \includegraphics[scale=0.45]{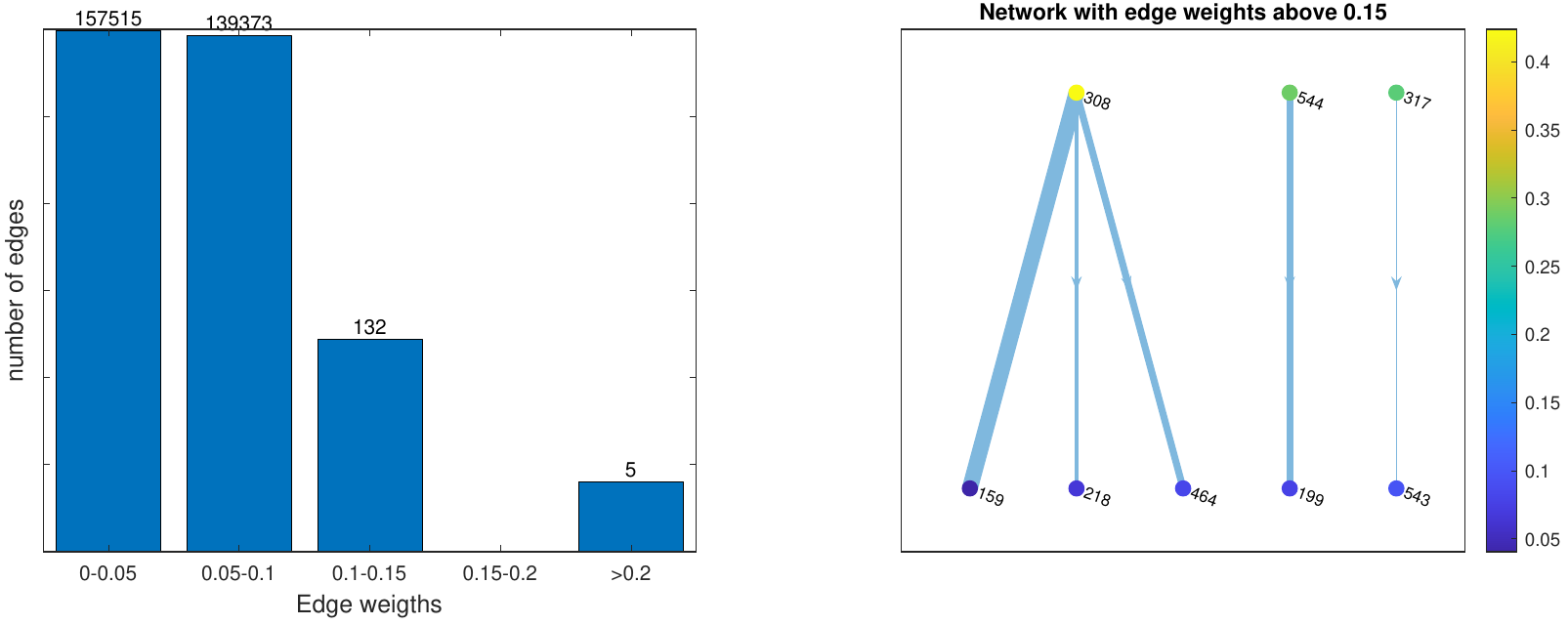}
    \caption{(Left) Histogram of the edge weights, i.e., $\alpha_{ij}$ for the large email network. (Right) The subnetwork with the edges weight above 0.15.}
    \label{fig:Email_influence_network}
\end{figure}

We test only ExPKF for this experiment since it requires less computer memory and runs faster than MM on our computational resources. We set the initial values $\alpha^{ij}=0.1$ and initialise $\mu^i$ by the average count on the $i-$th node.
We set the decay rate $\beta=0.15$ for all nodes; we tested a few other values, and the results are qualitatively the same. We present the estimated influence network in Figure~\ref{fig:Email_influence_network}. It is clear that the network is extremely sparse. We can identify only 5 edges that would suggest a strong influence. Although the number of emails sent by node 308 is close to the median value, we can identify its relatively higher influence on a few other nodes, all of which have a low number of sent-out emails.

\section{Conclusion}\label{s:con}

This work presents a significant development in foundations and methods for reconstructing influence networks from a time-series of count data through parameter estimation of discrete-time, multivariate Hawkes or Cox processes. Developing methods for inference for count data is important as it is very common in applications when timestamp data is not available or does not make sense to collect, e.g. in epidemiology applications, but this area is significantly less developed than for timestamp data where, to the best of the authors' knowledge, there were previously no methods for dealing with count data. Despite count data having less information than the time-stamp data, we find that network reconstruction is still possible. We demonstrate an application of the ensemble-based EM algorithm for certain doubly-stochastic processes (such as Log-Gaussian Cox process) that can be presented in state-space form.
Our implementation is based on the forward filtering-backward smoothing procedure using the bootstrap particle filter for the forward filtering, which is followed by backward smoothing simulation. We demonstrated the that the Ensemble-EM method is able to carry out the network reconstruction through synthetic experiments with known ground truths for small networks. 
\par
We observe that the parameter settings in our synthetic experiments are carefully adjusted to ensure bounded intensity. When the intensity grows unbounded, network recovery becomes infeasible within our experimental framework. This challenge translates to practical scenarios where intensity values become excessively large. A similar issue is discussed in~\cite{Diouane, Rizoiu22}, where the process eigenvalue approaches 1. It is evident that further research advancements are necessary to address this limitation effectively.
\par
This paper lays the foundations for other smoothing methods that could be used instead of the forward-filtering-backward smoothing approach for the ensemble-based EM depending on the structure of the state model and the observational likelihood. For example, it was demonstrated in \cite{Moral2010ForwardSU} that it is possible to bypass entirely the backward smoothing to compute expectations in the setting of an online EM method. This would significantly reduce the memory storage requirements for the backward smoothing simulation and allow for larger networks to be handled. Future work will look at the development of the ensemble-based approximate filtering using similar concepts from the ensemble-based Kalman smoother (EnKS) developed in geophysical applications \cite{Evensen2018AnalysisOI}. The EnKS uses the ensemble to approximate the density of the one-step push-forward state. This ensemble is then updated to fit the observation problem under the approximately linear model and Gaussian observational noise. In the current context, however, the observation equation can be nonlinear in the parameters and may not be close to Gaussian for a time-series with small counts. Further study in this direction to improve the E-step of the ensemble-based EM will open up applications to the influence network reconstruction problem with more complicated state-space models.
\par
We then presented the MM-based algorithm and the ExPKF algorithm to handle real-world applications when the linear Hawkes model is a reasonable assumption. The MM algorithm is designed to handle large batch data. We select a tight upper bound so that each parameter can be updated separately in a parallel manner. The ExPKF algorithm is a sequential approach that assumes a known decay rate for the Hawkes model. This key assumption enables the rank-1 update in the algorithm to avoid the costly inversion of a large matrix and by estimating each node independently, our algorithm can be efficiently parallelized. Investigation of the ExPKF on synthetic data again showed excellent results in determining the hidden network structure. 
\par
We demonstrated the performance of the methods using numerical experiments with known ground truths for both perfect and imperfect model scenarios, and both ExPKF and MM algorithms can recover the influence network structure when compared with the ground truth with good estimates of the strengths of the connections in the network. Several exciting areas for future research include looking at when the ExPKF algorithm becomes expensive for general Hawkes models where the inversion of the Hessian term in ExPKF cannot be performed via a rank-1 update; hence it can become a numerical issue for large-scale problems. For MM algorithms, a tight upper bound must be specifically designed for a given model. Therefore, for more general models, finding a tight upper bound allowing for parallel update of parameters is an interesting area to investigate. One thing this work opens up is the possibility of network reconstruction for applications in social networks and neural networks that hitherto remained out of reach with current methods. 

\appendix
\section*{Appendix A: Forward filtering-backward smoothing }
\par
We provide a brief review of the particle filtering (PF) and backward smoothing simulation (BSS) used to generate smoothed particles required to evaluate the surrogate function in \eqref{eq:Qfunsplit} for the ensemble-based EM algorithm. 
\par

\noindent\textbf{Particle Filter (PF):} Let $\xb_k^{f(\ell)}$ and $w_k^{f(\ell)}$ denote, respectively, the $\ell$-th (filtered) particle and its corresponding normalized weight at time step $k=0,1,\ldots, K$, where $K$ is the length of the time-series of count data.
\begin{enumerate}

\item Initialization: Randomly generate $N_f$ particles from an initial distribution
$\xb_0^{f(\ell)} \sim p\left(\xb_0\right)$ and set initial weights: $w_0^{f(\ell)}=1/N_f$ for $i=1,\ldots, N_f$.
\item Repeat this step for $k=1, \ldots, K$, 
\begin{enumerate}
    \item Draw random samples from the conditional predictive distribution, denoted by $\xb_k^{(\ell)}$ based on \eqref{eq:State}, i.e.,$\xb_k^{(\ell)}\sim N(\xb_k^{(\ell)};\Psi\left(\xb_{k-1}^{f(\ell)}\right),\Qb)$, and then generate the predictive conditional intensity $\lambda_k^{j,(\ell)}$ for all nodes $j=1,\ldots,m$ based on \eqref{eq:observationeq}
\item Update (unnormalized) weights based on the likelihood model
$$
\widetilde{w}_k^{f(\ell)} \propto{w}_{k-1}^{f(\ell)}\prod_{j=1}^m(\lambda_k^{j,(\ell)})^{\Delta N^i_k}\exp(-\lambda_k^{j,(\ell)}\delta t).
$$
and then normalize the weight by
${w}_k^{f(\ell)}:={w}_k^{f(\ell)}\left(\displaystyle\sum_{\ell=1}^{N}\widetilde{w}_k^{f(\ell)}\right)^{-1}.$
\item Perform resampling to add additional Monte Carlo variation when the effective sample size is low. We use the criteria below:
$$
N_{eff}:=\left(\displaystyle\sum_{\ell=1}^{N_f}\left(w_k^{(\ell)}\right)^2\right)^{-1}<0.5N_f.
$$
There are a number of methods for resampling. For simplicity, we use the systematic sampling algorithm described in \cite{Kitagawa1996MonteCF}, which costs $O(N_f)$. At the end of this step, we obtain $\xb_k^{f(\ell)}$ and $w_k^{f(\ell)}$.

\end{enumerate}
\end{enumerate}

\noindent\textbf{Backward Smoothing Simulation (BSS):}  

Let $\xb_k^{s(\ell)}$ denote the particle of the $\ell$-th smoothing path at time step $k=0$, $1, \ldots,K$.
\begin{enumerate}
    \item Initialization: Suppose $\xb_{0:K}^{f(\ell)}$ and $w_{0:K}^{f(\ell)}$, for $i=1, \ldots, N_f$, have been computed from (and stored during) the filtering process. Select $\xb_K^{s(\ell)}=\xb_K^{f(\ell)}$ with probability $w_K^{f(\ell)}$.
    \item Repeat this step (backward in time) for $k=K-1, \ldots, 0$,
        \begin{enumerate}
            \item For $\ell=1, \ldots, N_f$, calculate new weights according to~\eqref{eq:State}
                $$
                w_k^{s(\ell)} \propto w_k^{f(\ell)} p\left(\xb_{k+1}^{s(\ell)} \mid \xb_k^{f(\ell)}\right)=w_k^{f(\ell)}N\left(\xb_k^{f(\ell)};\xb_{k+1}^{s(\ell)},\Qb\right)
                $$
            \item Randomly select $\xb_k^{s(\ell)}=\xb_k^{f(\ell)}$ with probability $w_k^{s(\ell)}$.
            Repeat Step 1 and Step 2 $N_S$ times, where $N_s$ is the desired number of smoothing trajectories.
        \end{enumerate}

\end{enumerate}
The smoothing trajectories, $\xb_{0:K}^{s(\ell)}$ for $\ell=1,\ldots,N_s$ will then be used to estimate the parameters in the M-step of the EM algorithm, see again \eqref{eq:Qfunsplit}.
\section*{Appendix B: Maximization of $\mathcal{Q}_x$ in \eqref{eq:Qeval_LGCP}}
By neglecting $\mathcal{Q}_0$ in \eqref{eq:Qeval_LGCP}, we can find $\mu_1^{(\kappa+1)}, \omega_1^{(\kappa+1)}$ and $\epsilon^{(\kappa+1)}$ by maximizing $\mathcal{Q}_x$ only. Let $\beta = (1-\omega_1\delta t)$ and $\gamma=\omega_1\mu$. We can rewrite $\mathcal{Q}_x$ by
$$
\mathcal{Q}_x = -\frac{1}{2N_s\epsilon^2\delta t}\lVert \yb-\Ab\zb\rVert^2-\frac{1}{2}K\log\epsilon,
$$
where $\zb=\left[\beta,\gamma\right]^\top$, $\yb=\left[x_1^{s(1)},\ldots,x_K^{s(N_s)},\ldots,x_K^{s(1)},\ldots,x_K^{s(N_s)}\right]^\top$ and 
$$
\Ab^\top=\begin{pmatrix}
    x_0^{s(1)}&\cdots &x_{K-1}^{s(1)}&\cdots &x_0^{s(N_s)}&\cdots &x_{K-1}^{s(N_s)}\\
    \delta t&\cdots&\delta t&\cdots&\delta t&\cdots&\delta t 
\end{pmatrix}.
$$
Thus, maximizing $\mathcal{Q}_x$ is equivalent to finding $\zb$ to ``solves" the problem $\min\lVert \yb-\Ab\zb\rVert^2$, which is nothing but the normal equation if $\Ab^\top\Ab$ is full-rank or other techniques may be required to regularize the solution. However, since $\zb$ has to be positive, a quadratic programming should be used if unconstrained minimization fails to produce the desired positive solution.
\par
After obtaining $\zb^{(\kappa+1)}=\left[\beta^{(\kappa+1)},\gamma^{(\kappa+1)}\right]^\top$,  we can recover $\omega_1^{(\kappa+1)},\mu^{(\kappa+1)}$ from $\beta^{(\kappa+1)},\gamma^{(\kappa+1)}$. We also find the maximizing solution of $\epsilon$ by
$$
\epsilon^{(\kappa+1)} = \frac{1}{N_SK}\displaystyle\sum_{\ell=1}^{N_s}\displaystyle\sum_{k=1}^K\left(x_k^{s(\ell)}-\Psi_x(x_{k-1}^{s(\ell)})\right)^2,
$$
using $\omega_1^{(\kappa+1)},\mu^{(\kappa+1)}$ in $\Psi_x$ above.

\section*{Code availability}
Codes used to produce the results in this paper are available at:

\url{https://github.com/naratips/EM-MM-ExPKF}

\section*{Acknowledgement}
Naratip Santitissadeekorn is supported by the EPSRC grant EP/W02084X/1.

\bibliographystyle{plain}
\bibliography{ExPKFbib}

\end{document}